\pdfoutput=1
\RequirePackage{ifpdf}
\ifpdf 
\documentclass[pdftex]{sigma}
\else
\documentclass{sigma}
\fi

\numberwithin{equation}{section}

\newtheorem{Theorem}{Theorem}[section]
\newtheorem*{Theorem*}{Theorem}
\newtheorem{Corollary}[Theorem]{Corollary}
\newtheorem{Lemma}[Theorem]{Lemma}
\newtheorem{Proposition}[Theorem]{Proposition}
 { \theoremstyle{definition}
\newtheorem{Definition}[Theorem]{Definition}

\newtheorem{Remark}[Theorem]{Remark} }

\begin{document}
\allowdisplaybreaks

\newcommand{\arXivNumber}{2306.11707}

\renewcommand{\PaperNumber}{043}

\FirstPageHeading

\ShortArticleName{Hexagonal Circular 3-Webs with Reducible Polar Curves of Degree Three}

\ArticleName{Hexagonal Circular 3-Webs\\ with Reducible Polar Curves of Degree Three}

\Author{Sergey I. AGAFONOV}

\AuthorNameForHeading{S.I.~Agafonov}

\Address{Department of Mathematics, S\~ao Paulo State University-UNESP, S\~ao Jos\'e do Rio Preto, Brazil}
\Email{\href{mailto:sergey.agafonov@gmail.com}{sergey.agafonov@gmail.com}}

\ArticleDates{Received November 18, 2023, in final form June 04, 2025; Published online June 13, 2025}

\Abstract{The paper reports the progress with the classical problem, posed by Blaschke and Bol in 1938. We present new examples and new classifications of natural classes of hexagonal circular 3-webs. The main results is the classification of hexagonal circular 3-webs with reducible polar curves of degree 3 and description of hexagonal circular 3-webs admitting a one-parameter M\"obius symmetry.}

\Keywords{circular hexagonal 3-webs}

\Classification{53A60}

\section{Introduction}

The problem to describe hexagonal 3-webs formed by circles in the plane appeared in the first monograph on the web theory published by Blaschke and Bol in 1938 (see \cite[p.\ 31]{BB-38}). The authors presented an example with 3 elliptic pencils of circles, each pair of pencils sharing a~common vertex, and observed that one can construct hexagonal circular 3-webs from hexagonal linear 3-webs, completely described by Graf and Sauer \cite{GS-24} as being formed by tangents to a fixed curves of third class. The construction involves a central projection from a plane to a unit sphere followed by stereographic projection to a plane. The corresponding circular 3-webs were described earlier by Volk \cite{V-29} and Strubecker \cite{S-32}.

Stereographic projection puts the problem into a natural framework of the M\"obius geometry: instead of planar circular webs we study circular webs on the unit sphere, thus treating circles and straight lines on equal footing. M\"obius geometry assigns points outside the unit sphere to circles on this sphere: the assigned point is the polar point of the plane that cuts the circle on the sphere. This sphere is called also {\it Darboux quadric}. Thus any circular 3-web on the unit sphere determines locally 3 curve arcs outside the Darboux quadric, one arc per web foliation. Globally these arcs may belong to one irreducible algebraic curve. In what follows we call this set of polar points a {\it polar curve} of the web. For example, the polar curve of the hexagonal circular 3-web obtained from a linear 3-web is a planar cubic, possibly reducible. The polar curve of the cited example from \cite{BB-38} splits into 3 non-coplanar lines.

In the same year as the book \cite{BB-38} appeared, Wunderlich published a new remarkable example of hexagonal circular 3-web. Its polar curve splits into 3 conics lying in 3 different planes. Since through a point $p$ on the unit sphere pass the circles whose polar points are intersection of the web polar curve with the plane tangent to the unit sphere at $p$, the Wunderlich web is actually 6-web, containing 8 hexagonal 3-subwebs.

Wunderlich gave also a construction of hexagonal 3-webs whose polar curve splits into 3 non-coplanar lines, two being dual with respect to the Darboux quadric and the third joining them. These webs were later rediscovered by other authors.

Further, he presented the following way to construct hexagonal 3-webs: for any one-para\-met\-ric group acting in the plane, choose 2 transversally intersecting curve arcs that are also transversal to the group orbits; acting on the arcs by the group one gets 2 foliations; the third is composed by the group orbits. These 3 foliations compose a (local) hexagonal 3-web. Choosing a one-parameter group either of translations, or of dilatations, or of rotations and taking two intersecting circles (a straight line counts as a circle), we get circular hexagonal 3-webs.

 Blaschke was well aware of the difficulty of the posed problem and, in his last book on the web geometry \cite{B-55}, discussed the simpler problem of classifying hexagonal circular 3-webs whose polar curve splits into 3 non-coplanar lines. Note that to a line corresponds a pencil of circles that is {\it hyperbolic} if the line spears the Darboux quadric, {\it elliptic} if the line completely misses the Darboux quadric, or {\it parabolic} if the line touches the Darboux quadric.

By the year 1977, the list of 6 types (one from \cite{BB-38} and five indicated in \cite{W-38}) of circular hexagonal 3-webs whose polar curve splits into 3 non-coplanar lines was completed by Erdo\v{g}an~\cite{E-74} and Lazareva \cite{L-77}. The first attempt to prove that the list is actually complete was published in 1989 by Erdo\v{g}an \cite{E-89}. Based on direct computational approach, it did not provide the crucial computation: in fact, a modern computer systems for symbolic computations shows that there must be a mistake in the proof presented in \cite{E-89} (see concluding remarks for further detail).

The Erdo\v{g}an's claim was proved only in 2005 by Shelekhov \cite{S-05}. His insight was to look into the singular set of the webs: defined globally, the webs under study inevitably have singularities. Shelekhov considered the simplest possible singularities where two of the three circular foliations are tangent. It turns out that hexagonality imposes a strong restriction: locally, such singular set is either a circle arc of the 3rd foliation or the common circle arc of the first two. The restriction was rigid enough to obtain all the types on the list.

Five new types of hexagonal circular 3-webs were presented by Nilov in 2014 \cite{N-14}. Polar curves for four of them split into a line and a conic. The fifth example may be viewed as a 5-web whose polar curve is a union of a line and two conics. Taking the line and two arcs on different conics as the polar curve, one gets a hexagonal 3-subwebs.

One can not help to observe that the polar curves of all the known examples are algebraic. Motivated also by the dual reformulation of the Graf and Sauer theorem, we consider the following natural class of 3-webs: hexagonal circular 3-webs with polar curve of degree three. The main result of the paper is the complete classification of such webs with reducible polar curves.\looseness=-1

The case of planar polar curve follows immediately from the Graf and Sauer theorem: 3~points on the polar curve corresponding to 3~circles through a point~$p$ on the sphere are the ones where the plane, tangent to the sphere at $p$ meets the polar curve. This plane cuts the polar curve plane along the line. On the polar curve plane we get the configuration dual to the Graf and Sauer theorem.

The case of non-planar set of 3 lines was finally settled by Shelekhov \cite{S-05}.

We obtain a classification of 3-webs whose planar polar curve splits into a line and a smooth conic. Up to M\"obius transformation, there are 15 types, most of them depending on one parameter. Four types of five in Nilov's paper \cite{N-14} are webs of this list, namely, of the types 6, 10, 11 and 15, presented in Section \ref{split}. (In fact, Nilov has found only one M\"obius orbit from one-parametric family of orbits of our type 6.)

Another natural class that we study in this paper is the set of hexagonal circular 3-webs symmetric by action of one-parameter subgroup of the M\"obius group. We also give a complete classification of such webs.

To select candidates for hexagonal webs we exploit further the above mentioned observation of Shelekhov on simplest singularities of hexagonal 3-webs. The proof of the observed property in \cite{S-05}, based on considering the normal form of the web function is not complete: this normal form often does not exists at singular points (see concluding remarks for more detail). We make precise the ideas about the type of singularities and then prove the key singularity property.

For completeness, we also present the classification of hexagonal webs with 3 non-coplanar polar lines. The proof mainly follows the line taken by Shelekhov in \cite{S-05}.

\section{Hexagonal 3-webs, Blaschke curvature, singularities}
A planar 3-web ${\cal W}_3$ in a planar domain is a superposition of 3 foliations $\mathcal{F}_i$, which may be given by integral curves of three ODEs
$
\sigma _1=0$, $\sigma _2=0$, $ \sigma _3=0$,
where $\sigma _i$ are differential one-forms. At non-singular points, where the kernels of these forms are pairwise transverse, we normalize the forms so that
$\sigma_1+\sigma_2+\sigma_3=0.$ The {\it connection form} of the web ${\cal W}_3$ is a one-form $\gamma$ determined by the conditions
${\rm d}\sigma_i+\gamma \wedge \sigma_i=0$, $i=1,2,3$.
The connection form depends on the normalization of the forms $\sigma_i$, the {\it Blaschke curvature} $d\gamma$ does not.
\begin{Definition}
A 3-web is hexagonal if for any non-singular point there are a neighbourhood and a local diffeomorphism sending the web leaves of this neighbourhood in 3 families of parallel line segments.
\end{Definition}
Topologically, hexagonality means the following incidence property that has given its name to the notion: for any point $m$, each sufficiently small curvilinear triangle with the vertex $m$ and sides formed by the web leaves, may be completed to the curvilinear hexagon, whose sides are web leaves and whose ``large'' diagonals are the web leaves meeting at $m$ (see the gallery of pictures illustrating hexagonal webs in the next section). Computationally, hexagonality amounts to vanishing of the Blaschke curvature \cite{BB-38}.

Up to a suitable affine transformation, the forms $\sigma_i$ may be normalized as follows:
\[
\sigma_1=(Q - R)({\rm d}y - P{\rm d}x), \qquad \sigma_2=(R - P)({\rm d}y - Q{\rm d}x), \qquad \sigma_3=(P - Q)({\rm d}y - R{\rm d}x),
\]
where $P(x,y)$, $Q(x,y)$, $R(x,y)$ are the slopes of the tangent lines to the web leaves at $(x,y)$. Vanishing of the curvature writes as
\begin{gather}
(R-Q )[P_{xx} +(Q + R)P_{xy} +QRP_{yy}]+(P-R)[Q_{xx} +(P + R)Q_{xy}+PRQ_{yy}]\nonumber\\
\qquad{}+ (Q - P)[R_{xx} +(P + Q)R_{xy}+PQR_{yy}]\nonumber\\
\qquad{}+ \frac{(Q - R)(2P - Q - R)\bigl[P^2_x+(Q+R)P_xP_y+QRP^2_y\bigr]}{(P-Q )(P - R)}\nonumber \\
\qquad{} +
 \frac{(R-P )(2Q -P- R)\bigl[Q^2_x+(P+R)Q_xQ_y+PRQ_y^2\bigr]}{(Q - R)(Q-P)}
\nonumber \\
\qquad{}+ \frac{(P - Q)(2R-P-Q)\bigl[R_x^2+(P+Q)R_xR_y+PQR_y^2\bigr]}{(R-Q)(R - P)} \nonumber\\
\qquad{}+
\frac{(2R-P-Q)P_xQ_x }{P-Q}+
 \frac{(2P-Q-R)Q_xR_x}{Q - R}+
 \frac{(2Q-R-P)R_xP_x}{R-P} \nonumber\\
\qquad{}+ \frac{ \bigl(R^2-PQ\bigr)[P_xQ_y+P_yQ_x]}{(P - Q)} +
 \frac{\bigl(P^2 - QR\bigr)[Q_xR_y+Q_yR_x]}{Q-R} \nonumber\\
\qquad{} +
 \frac{\bigl(Q^2-PR\bigr)[R_xP_y+R_yP_x]}{R-P} + \frac{\bigl(2PQR -(P+Q)R^2\bigr)P_yQ_y}{Q-P}\nonumber\\
 \qquad{} +
\frac{\bigl(2PQR-(Q+R)P^2\bigr)Q_yR_y}{R-Q} +
\frac{\bigl(2PQR-(P+R)Q^2\bigr)R_yP_y }{P-R} =0.\label{curvature}
\end{gather}

If only one slope, say $R$, is given as an explicit functions of $x$, $y$ and $P$, $Q$ are roots of a~quadratic equation $P^2+AP+B=0$, then one finds the first derivatives of $P$, $Q$ by differentiating the Vieta relations
\begin{equation}\label{Vieta}
P+Q=-A,\qquad PQ=B,
\end{equation}
as a functions of $P$, $Q$ and the first derivatives of $A$, $B$. Differentiating these expressions, one gets also the second derivatives. Finally, excluding $P$ and $Q$ with the help of \eqref{Vieta}, one can rewrite \eqref{curvature} in terms of $A$, $B$, $R$ and their derivatives. The result is presented in Appendix~\ref{AppendixA}.

The webs considered later will inevitably have singularities: some kernels of the forms $\sigma_i$ can be not transverse or the forms can vanish at some points. We call a singular 3-web hexagonal if its Blaschke curvature vanishes identically at regular points. The simplest type of singularities of hexagonal 3-webs have the following remarkable property, first observed by Shelekhov \cite{S-05}.

\begin{Lemma}\label{sing}
Suppose that a hexagonal $3$-web, defined by three analytic direction fields $\xi_1$,~$\xi_2$,~$\xi_3$, has a singular point $p_0$ such that
\begin{enumerate}\itemsep=0pt
\item[$(1)$] all $\xi_i$ are well defined at $p_0$,
\item[$(2)$] $\xi_1$ and $\xi_2$ are transverse at $p_0$,
\item[$(3)$] $\xi_1=\xi_3$ at $p_0$,
\end{enumerate}
then either the leaves of $\xi_1$ and $\xi_3$ through $p_0$ coincide or
$\xi_1=\xi_3$ along the leaf of $\xi_2$ through $p_0$.
\end{Lemma}
\begin{proof}
The property is a consequence of separation of variables for hexagonal webs.
The second condition implies that we can rectify $\xi_1$ and $\xi_2$, i.e., choose some local coordinates $u$, $v$ so that $\xi_1=\partial_v$ and $\xi_2=\partial_u$. Then $\xi_3=-f(u,v)\partial_u+\partial_v$ with $f(u_0,v_0)=0$, where $p_0=(u_0,v_0).$
Consider the analytic functions $\varphi(u):=f(u,v_0)$.

If $\varphi(u)\equiv 0$, then $f(u,v_0)=0$ and the leaves of $\xi_1$ and $\xi_3$ are tangent along the line $v=v_0$, which is the leaf of $\xi_2$.

 If $\varphi(u)\not\equiv 0$, then $\varphi(u)=(u-u_0)^n\psi(u)$ with natural $n$ and $\psi(u_0)\ne 0$. Thus for any~$u_1$ close to $u_0$ with $u_1\ne u_0$ holds $f(u_1,v_0)=(u_1-u_0)^n\psi(u_1)\ne 0$. Now the hexagonality amounts to~${\partial_v\bigl(\frac{\partial_uf}{f}\bigr)=0}$ hence the germ $\tilde{f}$ of $f$ at $(u_1,v_0)$ factors $\tilde{f}(u,v)=a(u)b(v)$ with analytic germs $a$, $b$ at $u_1$ and $v_0$ respectively. One can choose $b(u)$ so that $b(v_0)=1$. Then $a(u)=\psi(u)=(u-u_0)^n\psi(u)$ is analytic also at $u_0$, the function of two variables $a(u)b(v)$ is analytic at $(u_0,v_0)$ and coincides with $f(u,v)$ at some neighborhood of $(u_1,v_0)$ included in the domain of $f$. Then by uniqueness $f(u,v)=a(u)b(v)$. Now observe that $a(u_0)=0$ which
implies $f(u_0,v)\equiv 0$ and the integral curves of $\xi_1$ and $\xi_3$ passing through $p_0$ coincide.
\end{proof}

\section{Projective model of M\"obius geometry}
Following Blaschke \cite{B-29}, we call the subgroup ${\rm PSO}(3,1)$ of projective transformations of $\mathbb{RP}^3$, leaving invariant the quadric
$
X^2+Y^2+Z^2-U^2=0$,
the M\"obius group.
For the reference, we present here infinitesimal generators of M\"obius group in homogeneous coordinates $[X:Y:Z:U]$ in $\mathbb{P}^3$, affine coordinates $x=\frac{X}{U}$, $y=\frac{Y}{U}$, $ z=\frac{Z}{U}$ in $\mathbb{R}^3$ and cartesian coordinates $(\bar{x},\bar{y})$ in $\mathbb{R}^2$ related to points $(x,y,z)$ on the unit sphere via stereographic projection
\smash{$
x=\frac{2\bar{x}}{1+\bar{x}^2+\bar{y}^2}$}, \smash{$ y=\frac{2\bar{y}}{1+\bar{x}^2+\bar{y}^2}$}, \smash{$z=\frac{1-\bar{x}^2-\bar{y}^2}{1+\bar{x}^2+\bar{y}^2}$}.
There are 3 rotations around the affine axes:
\begin{gather*}
R_z=Y\partial_X-X\partial_Y=y\partial_x-x\partial_y=\bar{y}\partial_{\bar{x}}-\bar{x}\partial_{\bar{y}},\\
R_x=Z\partial_Y-Y\partial_Z=z\partial_y-y\partial_z=\bar{x}\bar{y}\partial_{\bar{x}}+\frac{1}{2} \bigl(1-\bar{x}^2+\bar{y}^2\bigr)\partial_{\bar{y}},\\
R_y=X\partial_Z-Z\partial_X=x\partial_z-z\partial_x=-\frac{1}{2} \bigl(1+\bar{x}^2-\bar{y}^2\bigr)\partial_{\bar{x}}-\bar{x}\bar{y}\partial_{\bar{y}},
\end{gather*}
and 3 boosts (or ``hyperbolic rotations'')
\begin{gather*}
B_x=U\partial_X+X\partial_U=\partial_x-x(x\partial_x+y\partial_y+z\partial_z)=\frac{1}{2} \bigl(1-\bar{x}^2+\bar{y}^2\bigr)\partial_{\bar{x}}-\bar{x}\bar{y}\partial_{\bar{y}},\\
B_y=U\partial_Y+Y\partial_U=\partial_y-y(x\partial_x+y\partial_y+z\partial_z)= -\bar{x}\bar{y}\partial_{\bar{x}}+\frac{1}{2} \bigl(1+\bar{x}^2-\bar{y}^2\bigr)\partial_{\bar{y}},\\
B_z=U\partial_Z+Z\partial_U=\partial_z-z(x\partial_x+y\partial_y+z\partial_z)=-\bar{x}\partial_{\bar{x}}-\bar{y}\partial_{\bar{y}}.
\end{gather*}
The identity component of ${\rm PSO}(3,1)$ is well known to be isomorphic to the group $PSL_2(\mathbb{C})$, the isomorphism being given by the action $A(V)=AVA^*$ of $A\in SL_2(\mathbb{C})$ on the vector space of
matrices
\[
V=\left(
\begin{matrix}
X+U & Y+{\rm i}Z\\
Y-{\rm i}Z & U-X
\end{matrix}
\right)
\]
with real $X$, $Y$, $Z$, $U$. This action preserves determinant of~$V$, which is $U^2-X^2-Y^2-Z^2$.
By this isomorphism, the generators are represented by the following matrices:
\begin{gather*}
R_x=\frac{1}{2}\left(
\begin{matrix}
-{\rm i} & 0\\
0 & {\rm i}
\end{matrix}
 \right),\qquad
 R_y=\frac{1}{2}\left(
\begin{matrix}
0 & -{\rm i}\\
-{\rm i} & 0
\end{matrix}
 \right),\qquad
 R_z=\frac{1}{2}\left(
\begin{matrix}
0 & 1\\
-1 & 0
\end{matrix}
 \right),
\\
B_x={\rm i}R_x,\qquad B_y={\rm i}R_y, \qquad B_z={\rm i}R_z.
\end{gather*}

Two points $p_1=[X_1:Y_1:Z_1:U_1]$ and $p_2=[X_2:Y_2:Z_2:U_2]$ in $\mathbb{RP}^3$ determine a line with the Pl\"ucker coordinates
\begin{alignat*}{4}
& a:=X_1U_2-X_2U_1,\qquad&& b:=Y_1U_2-Y_2U_1, \qquad&& c:=Z_1U_2-Z_2U_1,&
\\
& f:=Y_1Z_2-Y_2Z_1,\qquad && g:=Z_1X_2-Z_2X_1, \qquad&& h:=X_1Y_2-X_2Y_1.&
\end{alignat*}
By direct computation, one proves the following fact.
\begin{Lemma}
All points of a line with Pl\"ucker coordinates $[a:b:c:f:g:h]$ are stable with respect to subgroup with the infinitesimal generator $aR_x+bR_y+cR_z+fB_x+gB_y+hB_z$.
\end{Lemma}
Observe that the line dual to $[a:b:c:f:g:h]$ is the one with coordinates $[-f:-g:-h:a:b:c]$, which corresponds to multiplication by ${\rm i}$ of the corresponding matrix representation of the generator.

A line in $\mathbb{RP}^3$ can be hyperbolic, elliptic or parabolic with the respect to the Darboux quadric. Considering simple representatives of these classes one sees that
\begin{itemize}\itemsep=0pt
\item[(1)] For hyperbolic line, the action of the corresponding operator on the dual line is M\"obius conjugate to the action of $R_z$ on the line $Z=U=0$ at infinity thus not having extra stable points in $\mathbb{RP}^3$.

\item[(2)] For elliptic line, the action of the corresponding operator on the dual line is M\"obius conjugate to the action of $B_z$ on the line $X=Y=0$ (which action leaves invariant two extra points $[0:0:\pm1:1]$).

\item[(3)] For parabolic line, the corresponding operator moves all points on the dual line (consider~${\partial_y=R_x+B_y}$).
\end{itemize}

In the following proposition we summarize further properties of the above correspondence.
\begin{Proposition}\label{operline}
Let $\xi =aR_x+bR_y+cR_z+fB_x+gB_y+hB_z$ be an infinitesimal operator of the M\"obius group.
\begin{enumerate}\itemsep=0pt
\item[$(1)$] The corresponding action of the one-parameter group is not loxodromic and therefore M\"obius equivalent $($i.e., conjugated$)$ either to rotation, or dilatation, or translation of the plane $(\bar{x},\bar{y})$ if and only if \begin{equation}\label{Pluck}
af+bg+ch=0.
\end{equation}
\item[$(2)$] This action is M\"obius equivalent to rotation if and only if \eqref{Pluck} and
$a^2+b^2+c^2>f^2+g^2+h^2$ are true, the line with Pl\"ucker coordinates $[a:b:c:f:g:h]$ being the set of polar points of the circular orbits.
\item[$(3)$] This action is M\"obius equivalent to dilatation if and only if~\eqref{Pluck} and
$a^2+b^2+c^2<f^2+g^2+h^2$ are true, the line with Pl\"ucker coordinates $[a:b:c:f:g:h]$ being the set of polar points of the orbits.
\item[$(4)$] The action is M\"obius equivalent to translation if and only if \eqref{Pluck} and
\begin{equation}\label{quadrcomplex}
a^2+b^2+c^2=f^2+g^2+h^2
\end{equation} are true, the line with Pl\"ucker coordinates $[a:b:c:f:g:h]$ being the set of polar points of the orbits.
\end{enumerate}
\end{Proposition}
\begin{proof}
The type of subgroup transformations is determined by the eigenvalues of the matrix representation for the generating operator $\xi$. Consider the characteristic polynomial for $\xi$
\[
\Psi (\lambda)=\lambda^2 +\frac{1}{4}\bigl( a^2 + b^2 + c^2 - f^2 - g^2 - h^2\bigr)+\frac{\rm i}{2}(af + bg +ch)
\]
and the matrices for non-loxodromic M\"obius representatives $R_z$, $B_z$ and $\partial_y=R_x+B_y$.
\end{proof}

\section{Polar curve splits into 3 non-coplanar lines}
Recall that {\it limit circles} of a hyperbolic pencil correspond to intersection points of the hyperbolic line with the Darboux quadric under the stereographic projection. {\it Vertexes} of an elliptic pencil are two points common to all the circles of the pencil, the vertexes correspond to the intersection of the Darboux quadric with the line, dual to the polar elliptic line. Considering a parabolic pencil as a limit case of elliptic, one calls the point corresponding to the point of tangency of the Darboux quadric with the parabolic line also the {\it vertex}. The vertex of a parabolic pencil is the common point for all circles of this parabolic pencil.

The reader may visualize the pencil fixed by a line $L$ thinking of its circles as cut on the Darboux quadric by the pencil of planes containing the dual line $L^*$.

Let us list hexagonal 3-webs with non-planar polar lines. There are 9 types of such webs up to
M\"obius transformation.
\begin{figure}[ht]
\includegraphics[width=0.35\textwidth, trim={4.5cm 4.5cm 4.5cm 4.5cm}, clip]{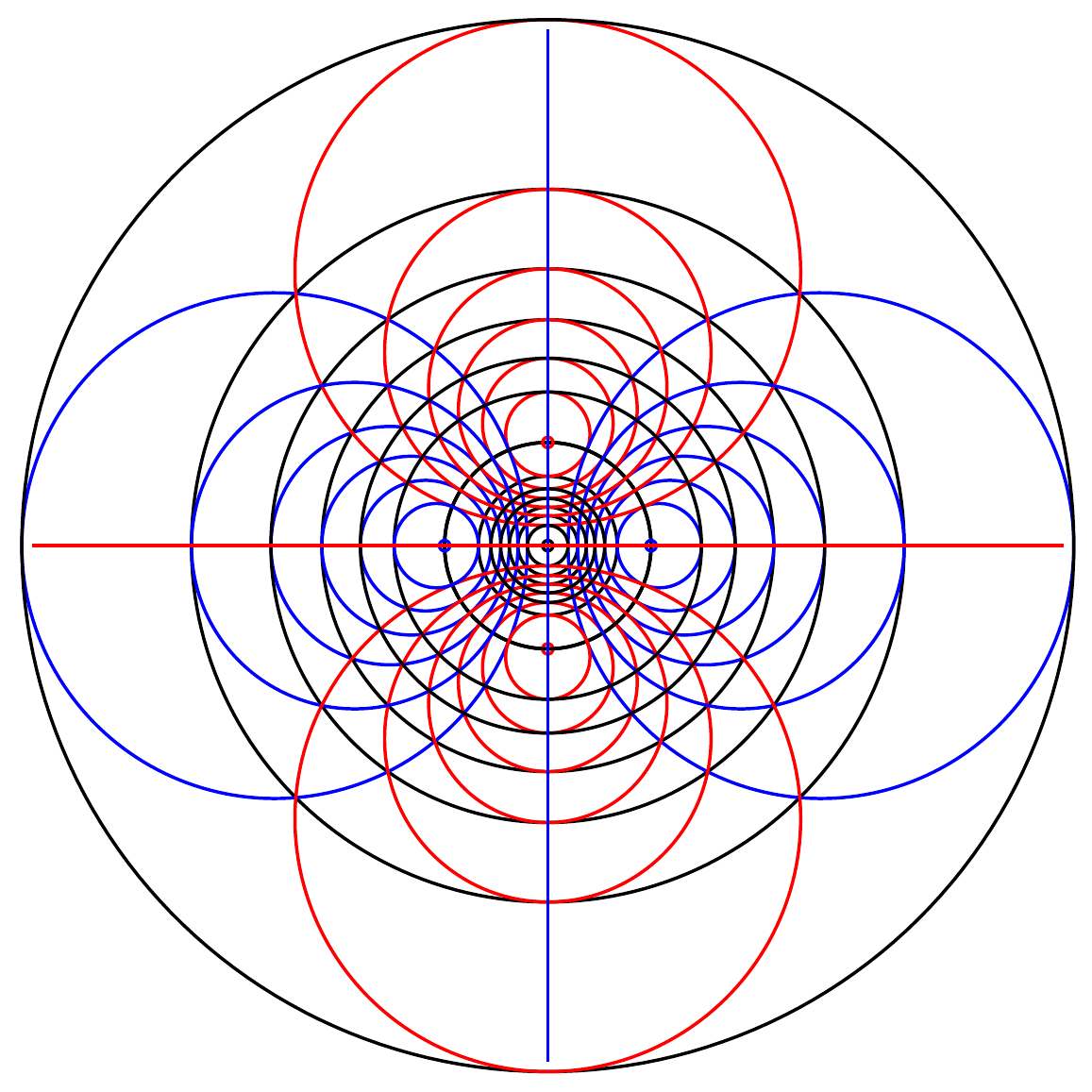}
\hspace{0.2cm}
\includegraphics[width=0.36\textwidth,trim={4cm 4cm 4cm 4cm}, clip]{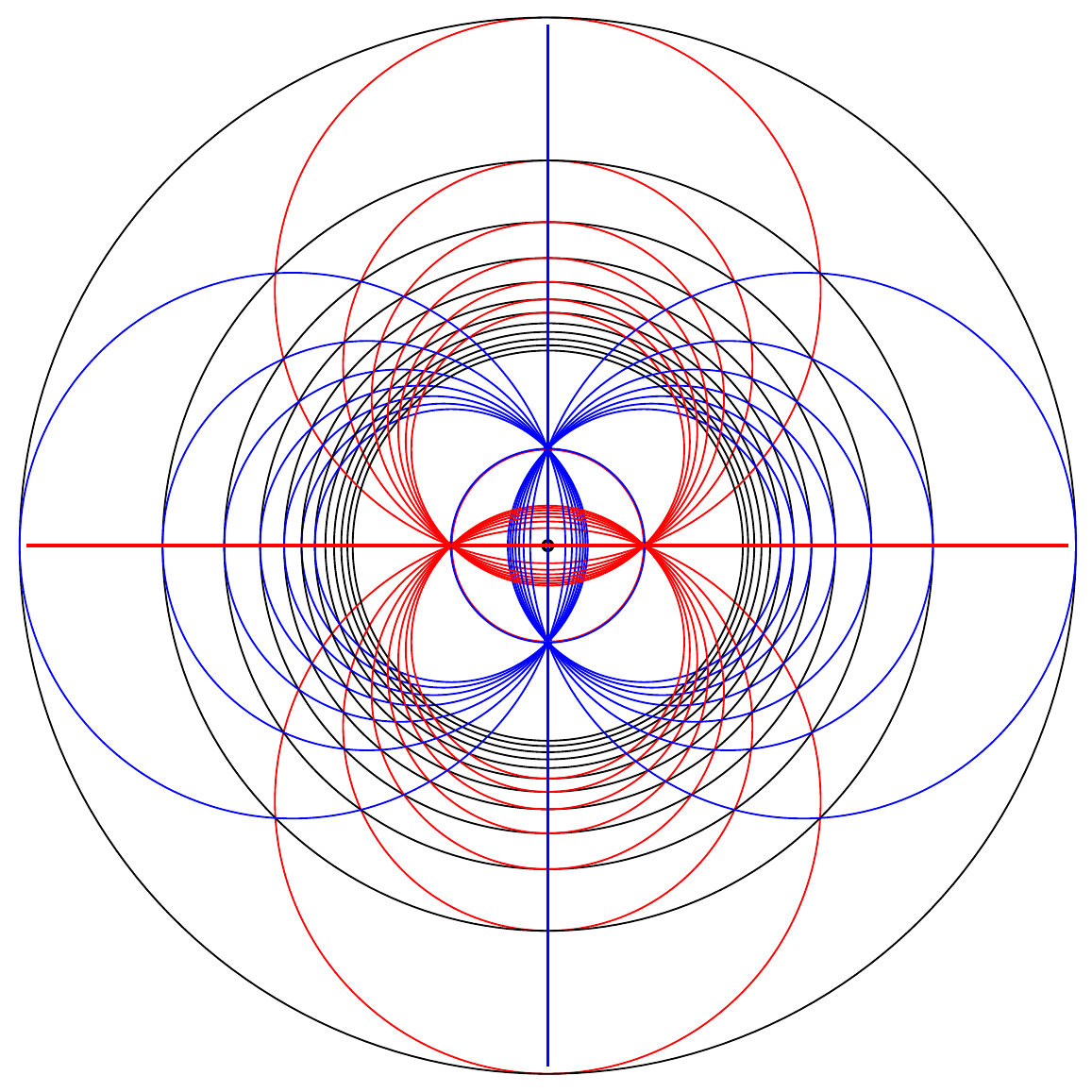}
\hspace{0.2cm}
\includegraphics[angle=90, width=0.24\textwidth, trim={1cm 4cm 1cm 4cm}, clip]{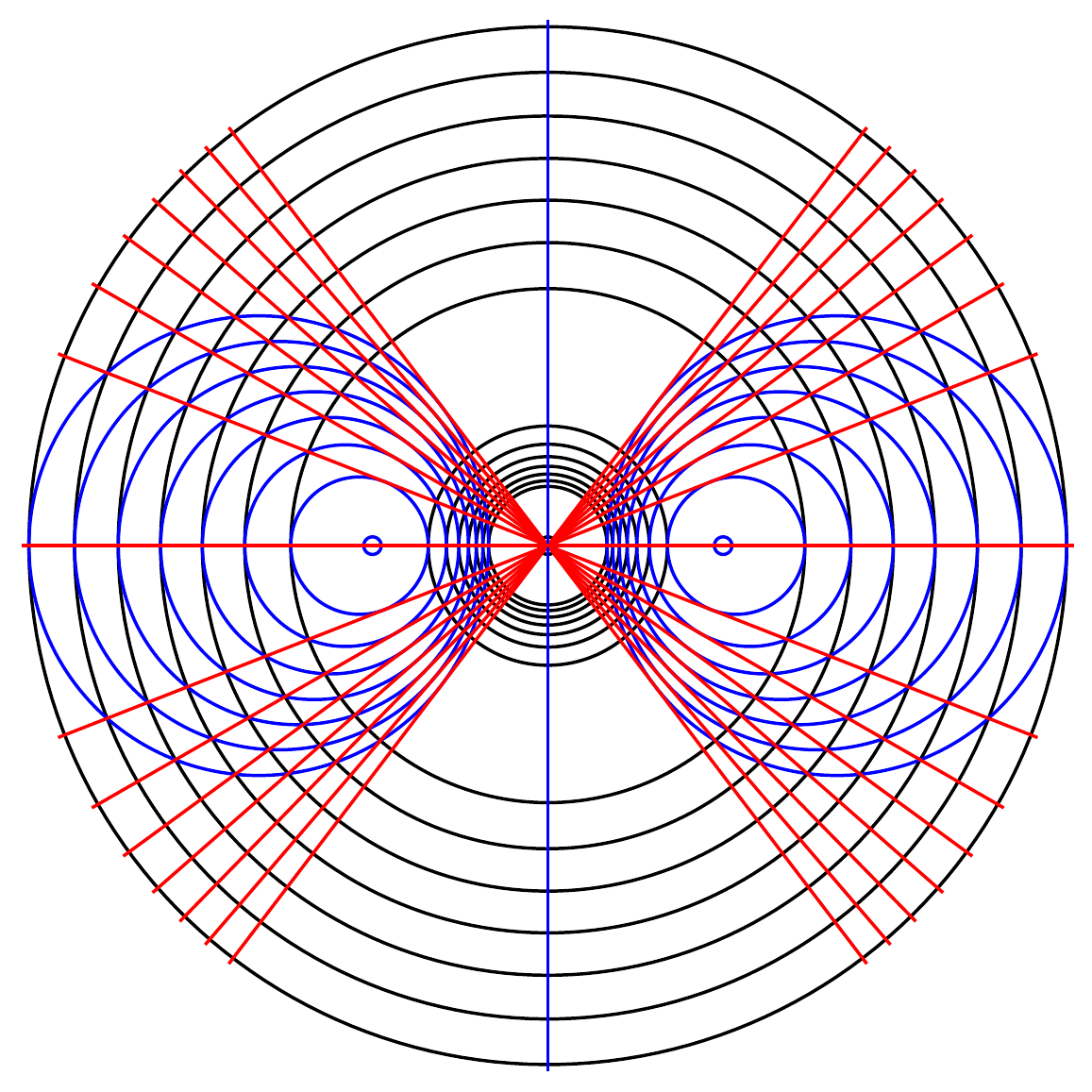}
 \caption{Hexagonal circular 3-webs. 3 hyperbolic pencils (left), 1 hyperbolic and 2 elliptic pencils (center), 1 elliptic and 2 hyperbolic pencils (right).}\label{f1}
\end{figure}

Take 3 polar lines intersecting inside the Darboux quadric so that each of the three lines contains the point dual to the plane of the other two polar lines, i.e., each pencil has a circle orthogonal to all the circles of the other two pencils. A representative of this web orbit, having one limit circle at infinity, is shown in Figure~\ref{f1} on the left. This type was described by Lazareva~\cite{L-77}.

Replacing two pencils by their orthogonal, we get the web in the center of Figure~\ref{f1}. In the projective model, we replace two polar lines by their dual ones.
 This web was also described by Lazareva~\cite{L-77}.

 \begin{figure}[ht]
\includegraphics[angle=90, width=0.280\textwidth, trim={3cm 4cm 3cm 4cm}, clip]{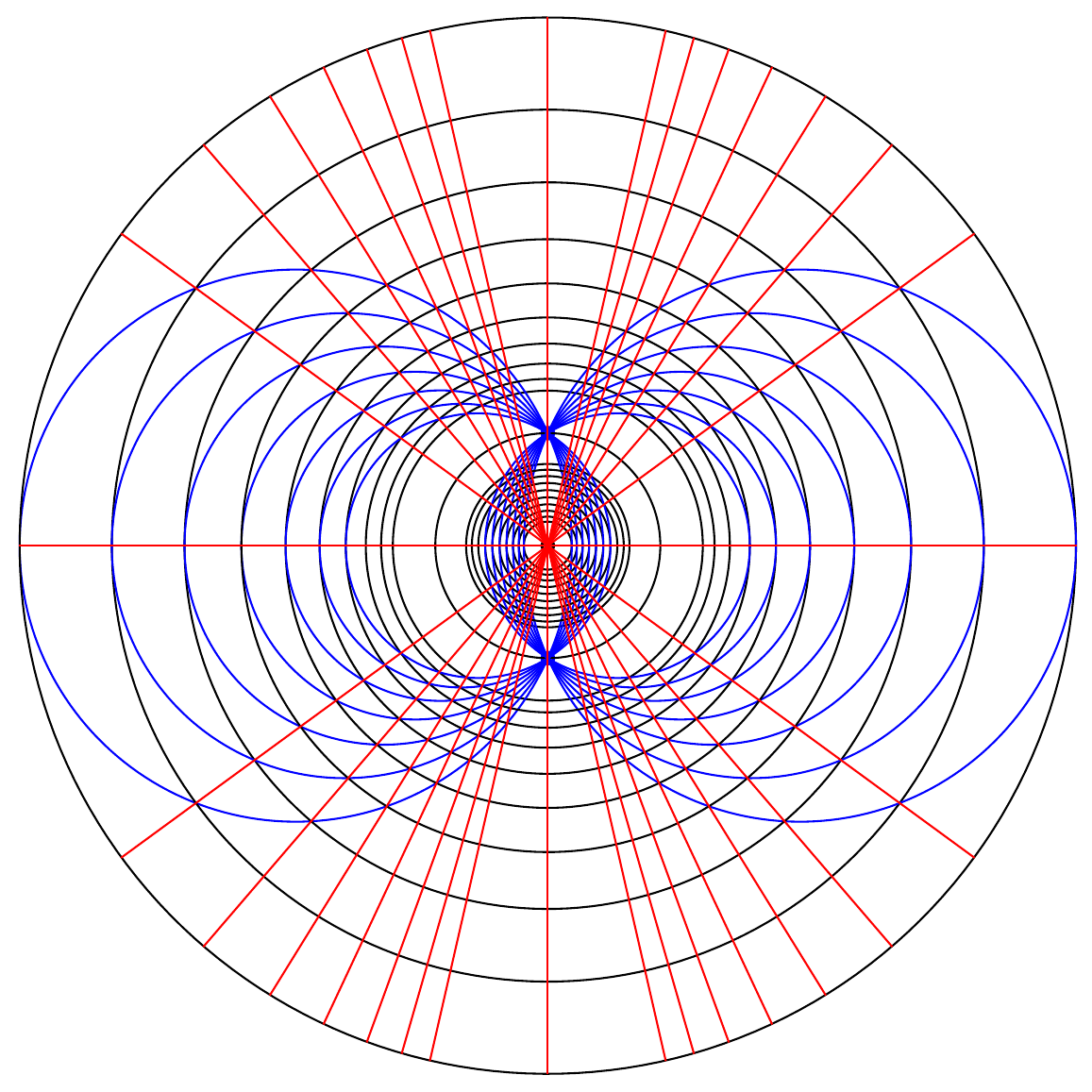}
\hspace{0.2cm}
\includegraphics[width=0.33\textwidth, trim={4cm 4cm 4cm 4cm}, clip]{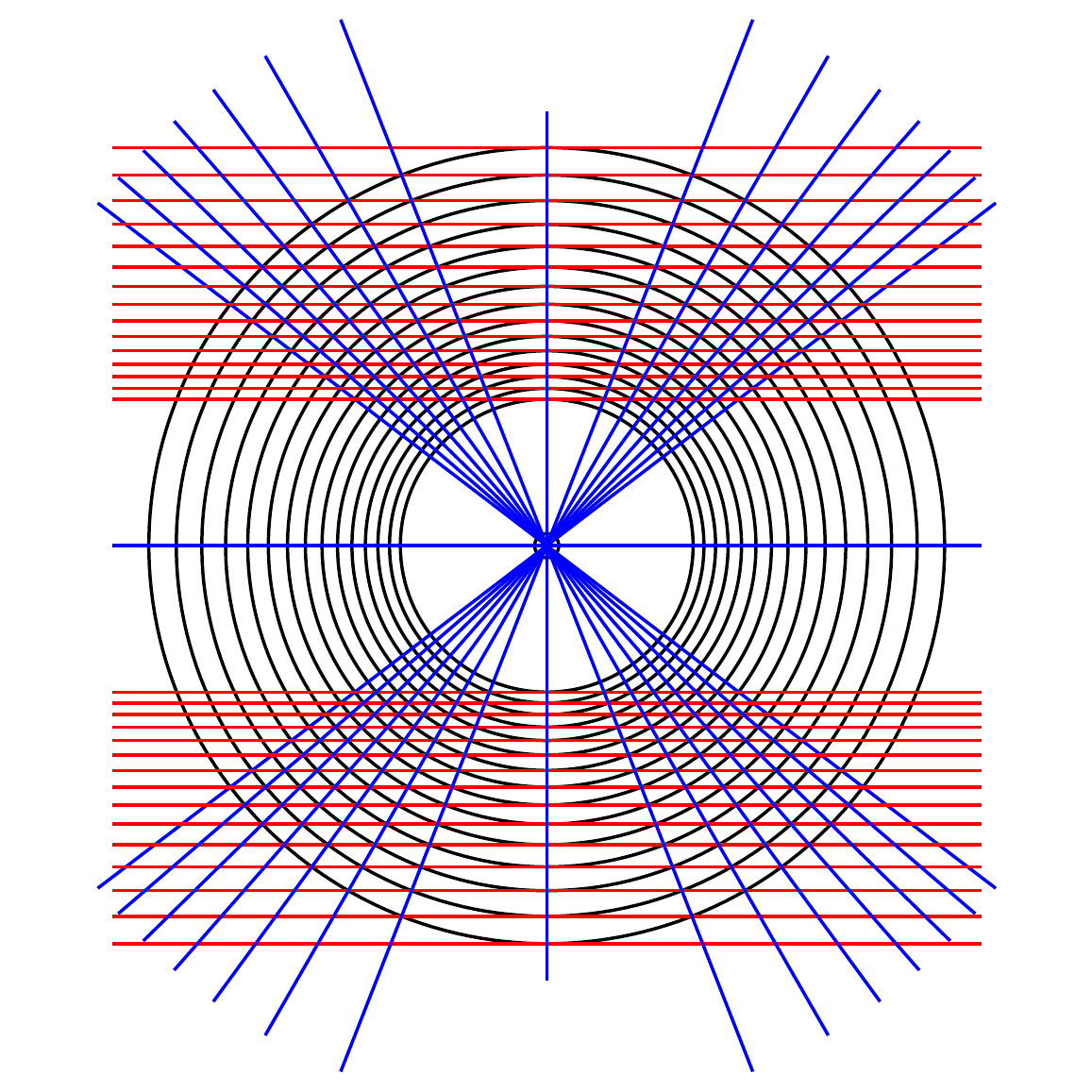}
\hspace{0.2cm}
\includegraphics[width=0.33\textwidth, trim={3.1cm 3.1cm 3.1cm 3.1cm}, clip]{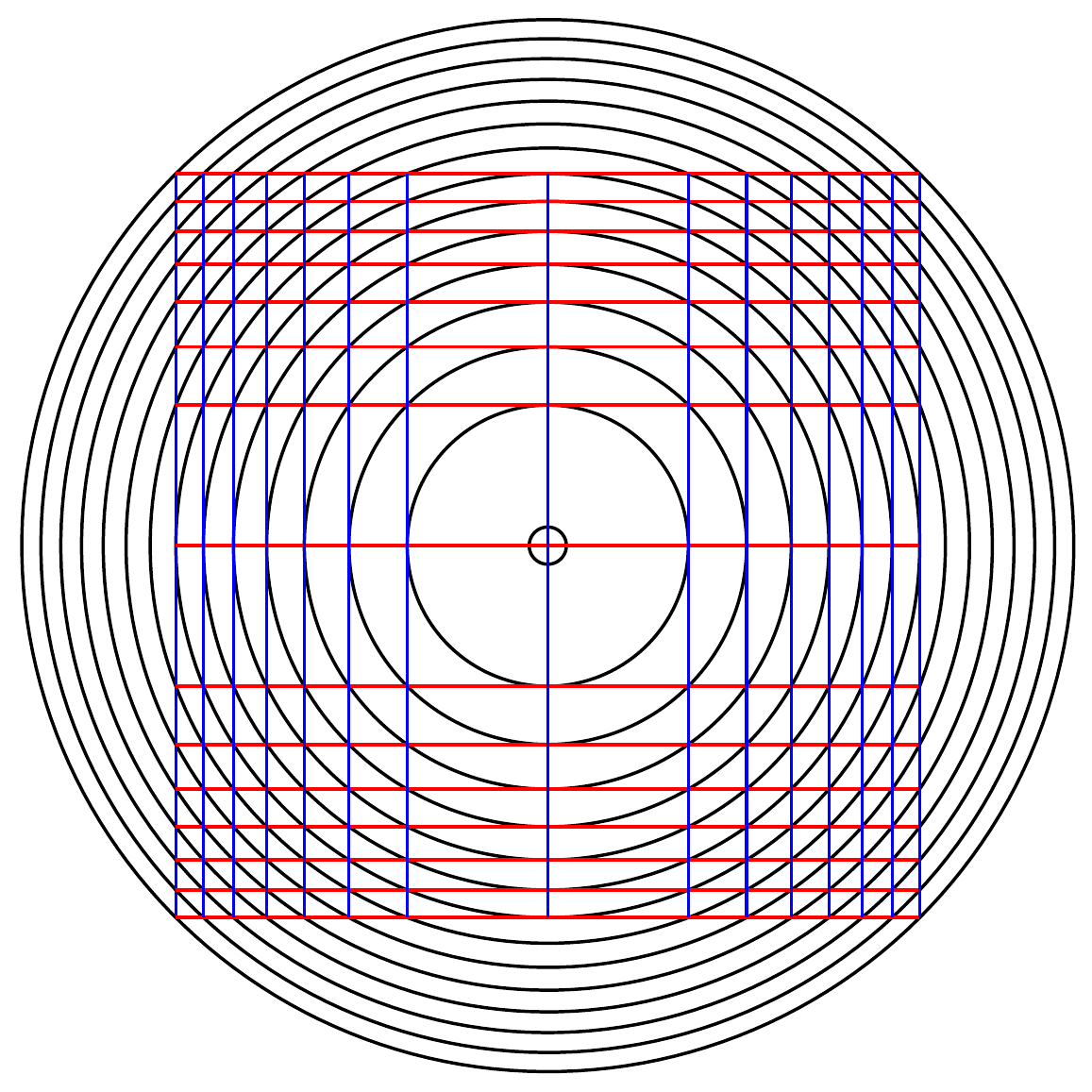}
 \caption{Hexagonal circular 3-webs. 1 hyperbolic and 2 elliptic pencils (left), 1 parabolic, 1 elliptic, and 1 hyperbolic pencil (center), 1 hyperbolic and 2 parabolic pencils (right).}\label{f2}
 \end{figure}

 Wunderlich \cite{W-38} mentioned the following construction, used later also by Balabanova and Erdo\v{g}an, to produce hexagonal 3-webs: take two dual polar lines and supplement it by a third intersecting that dual pair. There are four webs in the list, obtained in this way (see also \cite{B-73}): with two hyperbolic and one elliptic pencils on the right of Figure~\ref{f1};
with one hyperbolic and two elliptic pencils on the left of Figure~\ref{f2}; with one hyperbolic, one elliptic, and one parabolic pencils in the center of Figure~\ref{f2}; and with one hyperbolic and two parabolic pencils on the right of Figure~\ref{f2}.\looseness=-1

\begin{figure}[h]
\includegraphics[angle=90, width=0.31\textwidth, trim={8cm 8.2cm 8cm 8.2cm}, clip]{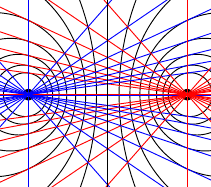}
\hspace{0.2cm}
\includegraphics[angle=90, width=0.28\textwidth, trim={5cm 6cm 5cm 6cm}, clip]{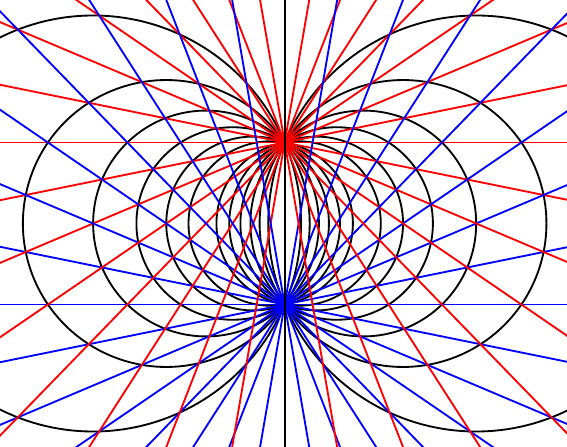}
\hspace{0.2cm}
\includegraphics[width=0.33\textwidth, trim={6cm 10.4cm 6cm 1.2cm}, clip]{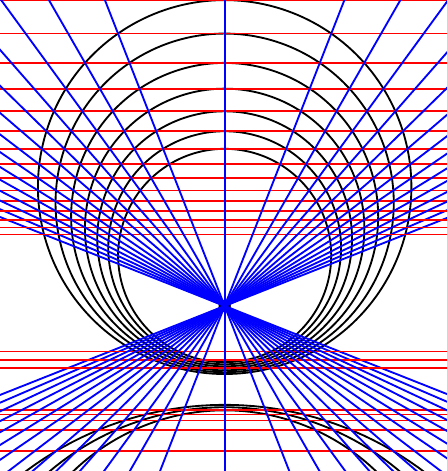}
 \caption{Hexagonal circular 3-webs. 1 hyperbolic and 2 elliptic pencils (left), 3 elliptic pencil (center), 1 hyperbolic, 1 elliptic and 1 parabolic pencil (right).}\label{f3}
\end{figure}

Another web with two elliptic and one hyperbolic pencils is depicted on the left of Figure~\ref{f3}. Its projective model has two elliptic polar lines, lying in a plane tangent to the Darboux quadric at a point, and a hyperbolic line passing through the points (different from the above tangency point), where duals to elliptic lines intersect the Darboux quadric. Its elliptic pencils share one common vertex at infinity, the other two vertices are also the limit circles of the hyperbolic pencil. This web was described by Erdo\v{g}an \cite{E-74}.

 In the center of Figure~\ref{f3}, there is a web with three elliptic pencils, the vertices of the pencils are two of three fixed points. It is historically the first hexagonal circular 3-web described in the literature \cite{BB-38}. The chosen representative has a vertex at infinity.

 Erdo\v{g}an \cite{E-74} found a web with one hyperbolic, one elliptic and one parabolic pencil, arranged so that the vertex $P$ of the parabolic pencil coincides with one vertex of the elliptic pencil, while the other vertex $E$ of the elliptic pencil coincides with one of the limiting circles of the hyperbolic, the common circle of the elliptic and the parabolic pencils being orthogonal to the circle passing through the second limiting circle of the hyperbolic pencil and the points $P$, $E$.
On the right of Figure~\ref{f3} is a M\"obius representative of this type web with $P$ at infinity. On the projective model, we have 3 pairwise distinct points on the Darboux quadric: $E$, $P$ and $H$. The hyperbolic polar line $L_h$ spears the quadric at $E$ and $H$, the parabolic polar $L_p$ touches the quadric at $P$ and intersects $L_h$, and the elliptic polar $L_e$ is dual to the line through $P$ and $E$ (and intersects $L_p$).\looseness=-1

\begin{figure}[h]\center
\includegraphics[width=0.3\textwidth, trim={6cm 0.4cm 6cm 12cm}, clip]{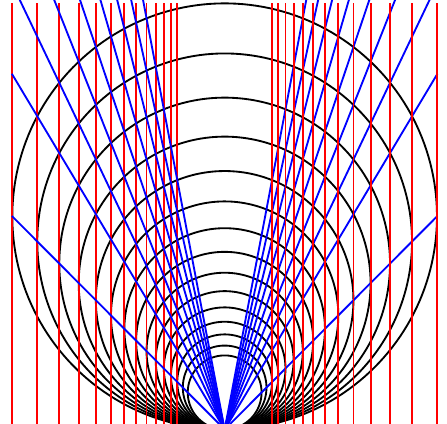}
 \caption{Hexagonal circular 3-web. 1 elliptic and 2 parabolic pencils.}\label{f4}
\end{figure}

Finally, we present a family of hexagonal 3-webs formed by 3 pencils, whose M\"obius orbits are parameterized by one parameter. Two pencils are parabolic with distinct vertexes, and the third is elliptic, whose dual to the polar line meets the Darboux quadric at these vertexes. A~representative of the family is shown in Figure~\ref{f4}, the vertexes being the origin and the infinite point. In this normalization one can fix the direction of one parabolic line, the direction of the other is arbitrary. Any web in this normalization is symmetric by dilatation $x\partial_x+y\partial_y$.

Surprising is not only the fact that the ``largest'' family was explicitly described only in 1977 by Lazareva \cite{L-77}, even more amazing is that the family falls within the general construction (see Introduction), which seems to have appeared first in the paper of Wunderlich \cite{W-38} in 1938!

In what follows we refer to the above described webs as {\it Blaschke--Wunderlich--Balabanova--Erdo\v{g}an--Lazareva list}.
To prove that there are no other classes, we will strongly use Lemma \ref{sing}. The singularity of the described type occurs when two of the three lines, tangent to Darboux quadric at a point and meeting each its own polar line, coincide and the point is not a vertex or limit circle of any pencil.
\begin{Proposition}\label{geometrypencil}
Consider the family of lines such that
\begin{itemize}\itemsep=0pt
\item[$(1)$] they meet two fixed lines $L_1$ and $L_2$, and

\item[$(2)$] they are tangent to the Darboux quadric.
\end{itemize}
 If the tangent points are on a circle, then either $L_1$ intersects $L_2$, or both $L_1$ and $L_2$ are tangent to the Darboux quadric, or one meets the dual of the other at some point on the Darboux quadric.

Moreover, in the case of skew $L_1$, $L_2$ tangent to the Darboux quadric at two points $p_1$, $p_2$, the curve of touching points splits into $2$ circles, their planes containing the line $p_1p_2$ and bissecting the angles between two planes $P_1$, $P_2$, where $P_i$ is the plane through $L_i$ and $p_1p_2$.
\end{Proposition}

\begin{proof}
Let $[a:b:c:f:g:h]$ be Pl\"ucker coordinates of a line $L$ touching the Darboux quadric. Then, by Proposition \ref{operline}, they satisfy equations \eqref{Pluck} and \eqref{quadrcomplex}. Therefore, $a^2+b^2+c^2\ne 0$ and we can normalize these coordinates to $a^2+b^2+c^2=f^2+g^2+h^2=1$. One easily calculates the point $p=(x,y,z)=(cg-bh,ah-cf,bf-ag)$ where the line $L$ touches the quadric. Due to normalization, one can rewrite this as
\begin{equation}\label{abcxyz}
a =hy -gz, \qquad b = f z - h x, \qquad c =gx -fy.
\end{equation}

 To simplify calculations, we can bring the Pl\"ucker coordinates of $L_1$ to simple form by M\"obius transformation. We suppose that $L_1$ and $L_2$ are skew since the case of $L_1$, $L_2$ intersecting is obvious.

If $L_1$ is hyperbolic, we can choose a representative as $L_1=[0:0:1:0:0:0]$. Then one has $m\ne 0$ since the lines $L_1$, $L_2=[u:v:w:k:l:m]$ do not intersect. (We used the fact that the Pl\"ucker coordinates of intersecting lines are orthogonal with respect to the bilinear symmetric form defining the Pl\"ucker quadric.) Let us set $m=1$. Moreover, the coordinates of $L_2$ satisfy the Pl\"ucker equation $ku+lv+w=0$ thus giving $w$. Since $L$ intersect $L_1$ and $L_2$, we have $h=0$ and $ka+lb+c+uf+gv+wh=0$, respectively. The above two equations cut a curve from the three-dimensional variety of lines touching the quadric.

The touching points $p$ trace a curve on the Darboux quadric.
This curve can be computed as follows. Equations \eqref{Pluck} and \eqref{abcxyz} imply $fx + gy + hz=0$, with $h=0$ we have ${g = -fx/y}$. Now normalization $f^2+g^2+h^2=1$ gives $\bigl(x^2 + y^2\bigr)f^2 = y^2$, which means $f$ is not identically zero. Therefore, intersection condition $ka+lb+c+uf+gv+wh=0$ is equivalent to $kxz + lyz + uy - vx - x^2 - y^2=0$. This equation cuts the curve of tangent points on the Darboux quadric (i.e., unit sphere centered at the origin).
If this curve is in a plane $z = Ax + By + C$, then by direct calculation one gets $C^2=1$. One can choose $C=1$ and then, by further calculations, we get~${A=-k}$,~${B=-l}$, $l=u$, $k=-v$. Thus $L_2=[u:v:0:-v:u:1]$ lies in the plane tangent to the Darboux quadric at $(0,0,-1)$, which is the point where the dual to $L_2$ spears the Darboux quadric. Note that the plane equation is $z = vx -uy + 1$.

If $L_1$ is hyperbolic, we choose $L_1=[0:0:0:0:0:1]$. Then $w\ne 0$ and we can normalize~${L_2=[u:v:1:k:l:m]}$ and the rest of the reasoning goes in a similar way.

If $L_1$ is parabolic, we choose $L_1=[0:-1:0:1:0:0]$. For skew $L_1$, $L_2$, we have $-l+u\ne 0$. Normalizing $l-u=1$ gives $l=u+1$. Further we rewrite \eqref{abcxyz} as
$
 f = b z - c y$, $ g = cx-a z$, $ h = a y - b x
 $
 and, proceeding as before, obtain by calculation that $L_2$ satisfy \eqref{Pluck}, which means it is tangent to the Darboux quadric.

 To check the last claim, we normalize the configuration so that $L_1$ and $L_2$ are tangent to the Darboux quadric at points $(0,0,\pm 1)$, look for a plane containing a circle of tangent points in the form $y=Ax$, find two solutions for $A$, and verify the geometry by calculation.
 \end{proof}

\begin{Definition}
 Motivated by Lemma \ref{sing}, we will refer to the circle of tangent points described in Proposition \ref{geometrypencil} whose polar point is not an intersection of $L_1$, $L_2$ as {\it singular circle}.
\end{Definition}

\begin{Remark}
 It is easy to see that the curve of touching points is of degree four. In the hypothesis of the proposition for the case of skew $L_1$, $L_2$ not tangent to the Darboux quadric, this curve also splits. One component is a circle and the other has just one real point, where dual to the elliptic line meets the hyperbolic.
 \end{Remark}

\begin{Corollary}\label{singcircle}
For hexagonal circular $3$-webs formed by three pencils with non-planar polar lines,
any two $L_1$, $L_2$ of polar lines are either dual or obey geometrical restriction described by Proposition {\rm\ref{geometrypencil}}. Moreover, the polar point of a singular circle, defined by two polar lines, belongs to the third one.
\end{Corollary}
\begin{proof}Let $L_1$, $L_2$ be skew but not dual and such that the curve of points, where lines $L$ meeting both $L_1$, $L_2$ touch the Darboux quadric, is not planar. Then by Lemma \ref{sing} any such line meets also the third polar line $L_3$. Since $L_1$, $L_2$ are skew, $L_3$ cannot intersect neither of~$L_1$,~$L_2$. Thus~$L_1$,~$L_2$,~$L_3$ belong to one ruling of some quadric $Q$ and the touching lines $L$ form the other ruling. But then $L_1$, $L_2$, $L_3$ are also tangent to the Darboux quadric. This contradicts our initial assumption.
\end{proof}

 Let $(p,q),(r,s)\in \mathbb{R}^2$ be vertices of an elliptic pencil, then the pencil circles form the family
\[
I(x,y):=\frac{[(p-x)(r-x)+(q-y)(s-y)]^2}{\bigl[(p-x)^2+(q-y)^2\bigr]\bigl[(r-x)^2+(s-y)^2\bigr]}=\text{const}.
\]
The circles are the integral curves of the ODE
\[
\omega_e:={\rm d}(I)=f(x,y){\rm d}x+g(x,y){\rm d}y=0.
\]
The circles of the hyperbolic pencil with limit circles at $(p,q)$, $(r,s)$ are orthogonal to the circles of the above elliptic one. They are the integral curves of the ODE
\[
\omega_h:=g(x,y){\rm d}x-f(x,y){\rm d}y=0.
\]

Now we work out the cases with parabolic pencils. Let $(p,q)\in \mathbb{R}^2$ be the vertex of a parabolic pencil and $[r:1-r]\in \mathbb{P}$ a direction orthogonal to the line tangent to all circles of the pencil. Then the pencil circles form the family
\[
\tilde{I}(x,y):=\frac{(x-p)^2 + (y-q)^2}{r(p-x) + (r - 1)(y - q)}=\text{const}.
\]
The circles are the integral curves of the ODE
$
\omega_p:=d\bigl(\tilde{I}\bigr)=0$.
For the exceptional direction~${(1,-1)}$, the pencil circles family is
\[
\bar{I}(x,y):=\frac{(x-p)^2 + (y-q)^2}{(x-p) - (y - q)}=\text{const}
\]
and the corresponding EDO
$
\omega_{\bar{p}}:={\rm d}\bigl(\bar{I}\bigr)=0$.

Observe that differential forms $\sigma_i$, describing a 3-web of circles formed by 3 pencils, are algebraic. Thus Lemma \ref{sing} remains valid also over complex numbers in passing from $\mathbb{RP}^3$ to~$\mathbb{CP}^3$. By circles here we understand sections of the complex Darboux quadric by complex planes. The complexification simplifies the proof of the following theorem.

\begin{Theorem}[\cite{S-05}]\label{3pencils}
Any hexagonal circular $3$-webs formed by three pencils with non-coplanar polar lines is M\"obius equivalent to one from the Blaschke--Wunderlich--Balabanova--Erdo\v{g}an--Lazareva list.
\end{Theorem}

\begin{proof}
We consider all types of non-coplanar polar line triples $L_1$, $L_2$, $L_3$.

$\bullet$ {\it Three hyperbolic pencils}. By Corollary~\ref{singcircle}, all lines intersect at one point $p$. This point cannot be outside the Darboux quadric. In fact, applying a suitable M\"obius transformation, we send this point to an infinite one, say $p_x=[1:0:0:0]$ and the plane of $L_1$, $L_2$ to the plane~${Z=0}$. Then, by Corollary~\ref{singcircle}, the polar line $L_3$ joins $p_x$ and $p_z=[0:0:1:0]$, the polar point of the plane $Z=0$. Thus $L_3$ cannot be hyperbolic as supposed. The point $p$ cannot be on the Darboux quadric either: we can send it to $p=(1,0,0)$ and the plane of $L_1$, $L_2$ to the plane $Z=0$. Now Corollary~\ref{singcircle} implies that $L_3$ joins $p$ and $p_z$. Thus $L_3$ is parabolic and not hyperbolic.

Therefore, $p$ is inside the Darboux quadric and one can send it to the origin $(0,0,0)$.
Corollary~\ref{singcircle} implies that any of the lines $L_1$, $L_2$, $L_3$ contains the point dual to the plane of the other two lines. Thus $L_1$, $L_2$, $L_3$ are orthogonal and
 we have the web shown on the left of Figure~\ref{f1}. Note that there is only one such web up to M\"obius transformation.

$\bullet$ {\it One elliptic and two hyperbolic pencils}.
We can suppose that $L_3$ is an elliptic polar line and that $L_3$ is the infinite line in the plane $Z=0$. Due to Corollary~\ref{singcircle}, the polar lines~$L_1$,~$L_2$, being hyperbolic, must intersect. Then the plane of $L_1$, $L_2$ is dual to some point on $L_3$ and therefore contains the line $X=Y=0$ dual to~$L_3$ and can be assumed to be the plane $Y=0$.

Suppose that none of $L_1$, $L_2$ is dual to $L_3$. If both $L_1$, $L_2$ intersect $L_3$, then the intersection point is $p_x=[1:0:0:0]$. Applying Corollary~\ref{singcircle} to the pair $L_1$, $L_3$, we conclude that $L_2$, joining $p_x$ and the polar point of the plane of $L_1$, $L_3$ does not meet the Darboux quadric and is not hyperbolic as supposed. Thus we can suppose that $L_1$, $L_3$ are skew and that $L_1$, by Corollary~\ref{singcircle}, contains $(0,0,-1)$. Since $L_1$ is not dual to $L_3$, we infer by Corollary~\ref{singcircle} that $L_2$ contains the dual point of the singular circle for $L_1$, $L_3$. This point lies in the tangent plane to the Darboux quadric at $(0,0,1)$, therefore $L_2$, being non-parabolic, cannot meet $L_3$ and cannot also contain $(0,0,1)$. Then it passes through $(0,0,-1)$, which contradicts the geometry restriction imposed by Corollary~\ref{singcircle}.

So we can assume that $L_2$ is dual to $L_3$, i.e., it is the line $X=Y=0$. Corollary~\ref{singcircle} prevents~$L_1$ to be skew with $L_3$. Thus it meets $L_3$ and therefore intersect $L_2$ in a point inside the Darboux quadric.
We obtain the hexagonal web equivalent to the one on the right of Figure~\ref{f1}.

$\bullet$ {\it One hyperbolic and two elliptic pencils}.
By Corollary~\ref{singcircle}, elliptic lines, say $L_2$, $L_3$, intersect. We can assume that the hyperbolic line $L_1$ is the coordinate axis $X=Y=0$.

First consider the case when two lines, $L_1$, $L_2$, are dual. Then $L_2$ is the infinite line in the plane $Z=0$ and we can assume the intersection line of $L_2$ and $L_3$ being the point $p_y=[0:1:0:0]$. Then $L_3$ cannot be skew with $L_1$ due to Corollary~\ref{singcircle}: the polar of the singular circle determined by $L_1$, $L_3$ is finite and cannot lie on the infinite line $L_2$. Thus $L_3$, being elliptic, intersect $L_1$ outside the Darboux quadric and we get the web equivalent to the one shown on the left of Figure~\ref{f2}.

Now suppose that no pair $L_1$, $L_i$ is skew. Then $L_2$ and $L_3$ intersect $L_1$. Since the triple $L_1$, $L_2$, $L_3$ is not coplanar, all three lines intersect at one point outside the Darboux quadric, which can be taken as $p_z=[0:0:1:0]$. By Corollary~\ref{singcircle}, the line $L_3$ contains the dual point of the singular circle determined by $L_1$, $L_2$ and the line $L_2$ contains the dual point of the singular circle determined by $L_1$, $L_3$. This fixes $L_1$, $L_2$, $L_3$ up to rotation around $z$-axis and gives the web shown in the center of Figure~\ref{f1}.

Finally, consider the case with skew but not dual $L_1$, $L_2$. Applying M\"obius transformation, we send the intersection point of $L_2$ and $L_3$ to $p_x=[1:0:0:0]$ (preserving the position of $L_1$). Then $L_3$ joins $p_x$ and the polar point of the singular circle determined by $L_1$, $L_2$. We obtain the web, equivalent to one on the left of Figure~\ref{f3}.

$\bullet$ {\it Three elliptic pencils}.
We treat this case using the complex version of Corollary~\ref{singcircle}. First, we conclude that, being non-coplanar, all three polar lines intersect at one point, which we can send to $p_y=[0:1:0:0]$. None of 3 real planes, containing a pair of polar lines, can miss completely the real Darboux quadric. In fact, if the real plane of $L_1$, $L_2$ do not intersect the real Darboux quadric then the polar point of the complex singular circle of $L_1$, $L_2$ is inside the real Darboux quadric and $L_3$, being elliptic, cannot contain this point. None of 3~real planes, containing a~pair of polar lines, can intersect the real Darboux quadric. If the real plane of~$L_1$,~$L_2$ cuts the real Darboux quadric, the polar point of the real singular circle of $L_1$, $L_2$ is outside the real Darboux quadric and the real plane of~$L_1$,~$L_3$ do not meet the real Darboux quadric. Thus all 3~planes are tangent to the Darboux quadric.
A representative of such web is shown in the center of Figure~\ref{f3}.

$\bullet$ {\it One parabolic and two hyperbolic pencils}.
By Corollary~\ref{singcircle}, all 3 polar lines must intersect in one point. The intersection point cannot be inside the Darboux quadric since one line is parabolic. It also cannot be outside: we can send it to $p_x=[1:0:0:0]$ and the parabolic line, joining $p_x$ and the point dual to the plane of hyperbolic lines, will miss the Darboux quadric. Therefore, this point is on the Darboux quadric. One can move the plane of hyperbolic lines to~${Z=0}$. Then the parabolic line contains $p_z=[0:0:1:0]$ and none of the hyperbolic lines can pass through the point dual the singular circle of the other two lines. Thus there is no hexagonal web with non-coplanar parabolic and two hyperbolic polar lines.

$\bullet$ {\it One parabolic, one hyperbolic and one elliptic pencil}.
By Corollary~\ref{singcircle}, the parabolic line meets the other two.
If the hyperbolic line intersects the elliptic at some point $p$ outside the Darboux quadric, we can move the hyperbolic line to $X=Z=0$ and $p$ to $p_y=[0:1:0:0]$. Then the point $p_s$, dual for the singular circle of hyperbolic and elliptic lines, is infinite and the third polar line, joining $p_y$ and $p_s$ cannot be parabolic.

Therefore, hyperbolic and elliptic lines are skew and Corollary~\ref{singcircle} fixes the configuration up to M\"obius transformation: if these lines are dual we get the type shown in the center of Figure~\ref{f2}, otherwise the type on the right of Figure~\ref{f3}.

$\bullet$ {\it Two parabolic and one hyperbolic pencils}.
By Corollary~\ref{singcircle}, the hyperbolic line meets both parabolic lines.

If the parabolic lines are skew, then the corresponding two singular circles have their polar points on the line dual to the one joining the points of tangency of parabolic lines with the Darboux quadric. This dual line is elliptic, therefore this configuration is not possible.

If the parabolic lines intersect outside the Darboux quadric, then we can bring the plane of their intersection to $Z=0$. Now the hyperbolic line contains $p_z=[0:0:1:0]$ by Corollary~\ref{singcircle} and, intersecting the both elliptic lines, must meet them at their common point. Then it misses the Darboux quadric and is not hyperbolic.

Therefore, the parabolic lines are tangent to the Darboux quadric at the same point. Since the polar lines are not coplanar the hyperbolic line contains this point. We can bring the hyperbolic line to $X=Y=0$. Then Corollary~\ref{singcircle} implies that the parabolic lines are dual and we obtain
the web type shown on the right of Figure~\ref{f2}.

$\bullet$ {\it Two parabolic and one elliptic pencils}.
By Corollary~\ref{singcircle}, the elliptic line meets both parabolic lines.

If the parabolic lines are skew, then the corresponding two singular circles have their dual points on the line dual to the one joining the points of tangency of parabolic lines with the Darboux quadric. The third polar line must contain these point, it is elliptic and we get
the web type presented in Figure~\ref{f4}.

If the parabolic lines intersect outside the Darboux quadric, then we can bring the plane of their intersection to $Z=0$. By Corollary~\ref{singcircle}, the elliptic line contains $p_z=[0:0:1:0]$ and, intersecting the both elliptic lines, must meet them at their common point. Thus we get the type shown in Figure~\ref{f4}.

If the parabolic lines are tangent to the Darboux quadric at the same point, then the third line, being elliptic, cannot pass through this point. Therefore, it is coplanar with the parabolic lines.

 $\bullet$ {\it Three parabolic pencils}.
Suppose that two polar lines $L_1$ and $L_2$ intersect. If the intersection point is outside the Darboux quadric, then we bring the plane of $L_1$, $L_2$ to $Z=0$. Since the third polar line $L_3$ does not lie in this plane it contains $p_z=[0:0:1:0]$ by Corollary~\ref{singcircle}. Therefore, it is skew with at least one of $L_1$, $L_2$. Let it be $L_1$. Then by Corollary~\ref{singcircle}, the line $L_2$ contains both dual point of two singular circles of $L_1$, $L_3$ which is obviously not possible.

If $L_1$, $L_2$ touch the Darboux quadric at the same point, then the line $L_3$ is skew with at least one of $L_1$, $L_2$. Again, this is precluded by Corollary~\ref{singcircle}.

Therefore, $L_1$, $L_2$, $L_3$ are pairwise skew. Consider two singular circles $C_1$, $C_2$ of $L_1$, $L_2$, their polar points $p_1$, $p_2$ and the family of lines $L$ touching the Darboux quadric and meeting both~$L_1$,~$L_2$. The third line $L_3$, being parabolic, cannot contain both points $p_1$, $p_2$. If $p_1 \notin L_3$, then by Lemma \ref{sing} the family of lines $L$ touching the Darboux quadric at points of $C_1$ must meet also~$L_3$. Therefore, the lines $L$ constitute one ruling of a quadric touching the Darboux quadric along~$C_1$. Therefore, $L_1$, $L_2$, $L_3$ belong to the second ruling. If $p_2\in L_3$, then, sending the plane of~$C_1$ to $Y=0$ and the points, where $L_1$, $L_2$ touch the Darboux quadric, to $(0,0,\pm 1)$, we see that~${p_2=[1:0:0:0]}$ since $C_2$ is orthogonal to $C_1$ and passes through $(0,0,\pm 1)$. Then~$L_3$ must touch the Darboux quadric also at one of the points $(0,0,\pm 1)$, which contradicts the initial assumption that all 3 lines are skew. Thus $p_2\notin L_3$. Then the family of lines $L$ touching the Darboux quadric at points of $C_2$ must meet also $L_3$. This is not possible as the lines $L$, constituting a ruling of a quadric that touch the Darboux quadric along $C_1$ cannot meet the orthogonal circle $C_2$.
\end{proof}

\section{Polar curve splits into conic and straight line}\label{split}

First we describe the types then we prove that the list is complete.
Some types are one-parametric families and we denote the parameter value by $c$. The polar conic will be given either by an explicit parametrization of the circle equations and we reserve $u$ for the parameter, or by indicating the conic equations. The former representation gives also a parametrization of the polar conic: to a circle
$
\epsilon \bigl(x^2+y^2\bigr)+\alpha x+\beta y+\gamma=0$
(where $\epsilon=0$ or $\epsilon=1$, the case $\epsilon=0$ giving a line)
corresponds the polar point with the tetracyclic coordinates
$
[\alpha:\beta:\gamma-\epsilon:-\gamma-\epsilon]$.
The parameter for circles in the pencil will be denoted by~$v$.

To check hexagonality, one computes the Blaschke curvature using the formula in Appendix~\ref{AppendixA} as follows. The polar conic gives a one-parameter family of circles on the unit sphere.
The stereographic projection from the ``south'' pole
\[
X=\frac{2x}{1+x^2+y^2}, \qquad Y=\frac{2y}{1+x^2+y^2}, \qquad Z=\frac{1-x^2-y^2}{1+x^2+y^2},
\]
 transforms this family into a family of circles in the plane parametrized by the points of the conic. The ODE for circles, defined by the conic,
\[
P^2+A(x,y)P+B(x,y)=0, \qquad P=\frac{{\rm d}y}{{\rm d}x},
\]
 is obtained by differentiating and excluding the coordinates of the conic points.
 The slope $R$ comes from the pencil of circles.

For some types, we add an additional geometric detail in the ``title'' to separate the types.

\begin{figure}[t]\centering
\includegraphics[width=0.30\textwidth, trim={6.1cm 2cm 0.2cm 2cm}, clip]{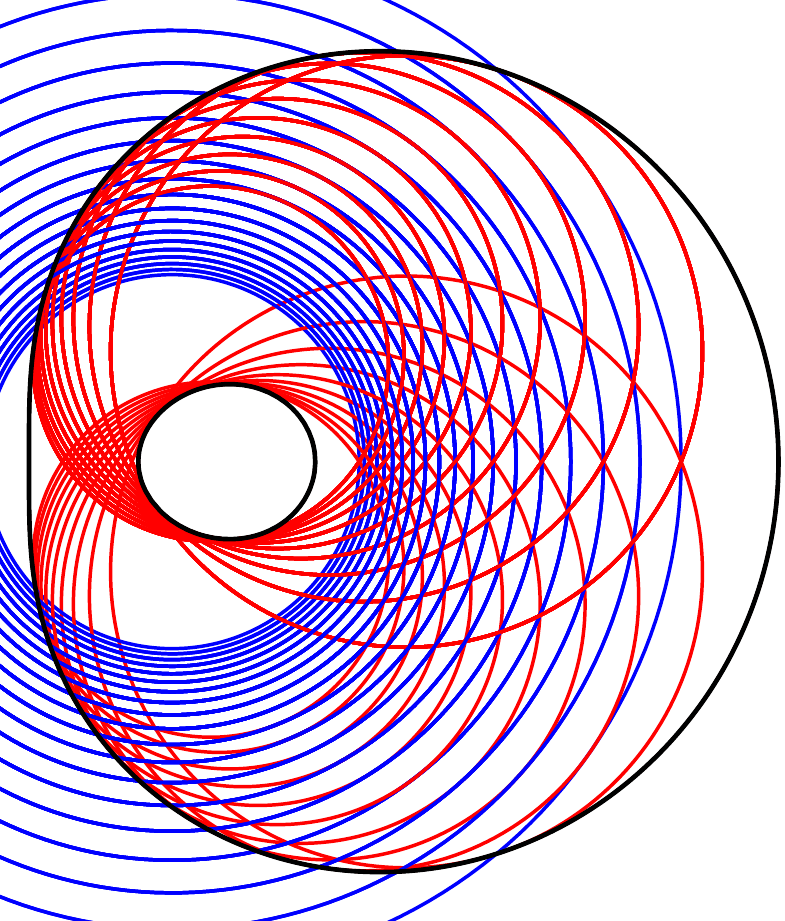}
\includegraphics[width=0.33\textwidth, trim={0.3cm 0.3cm 0.3cm 0.3cm}, clip]{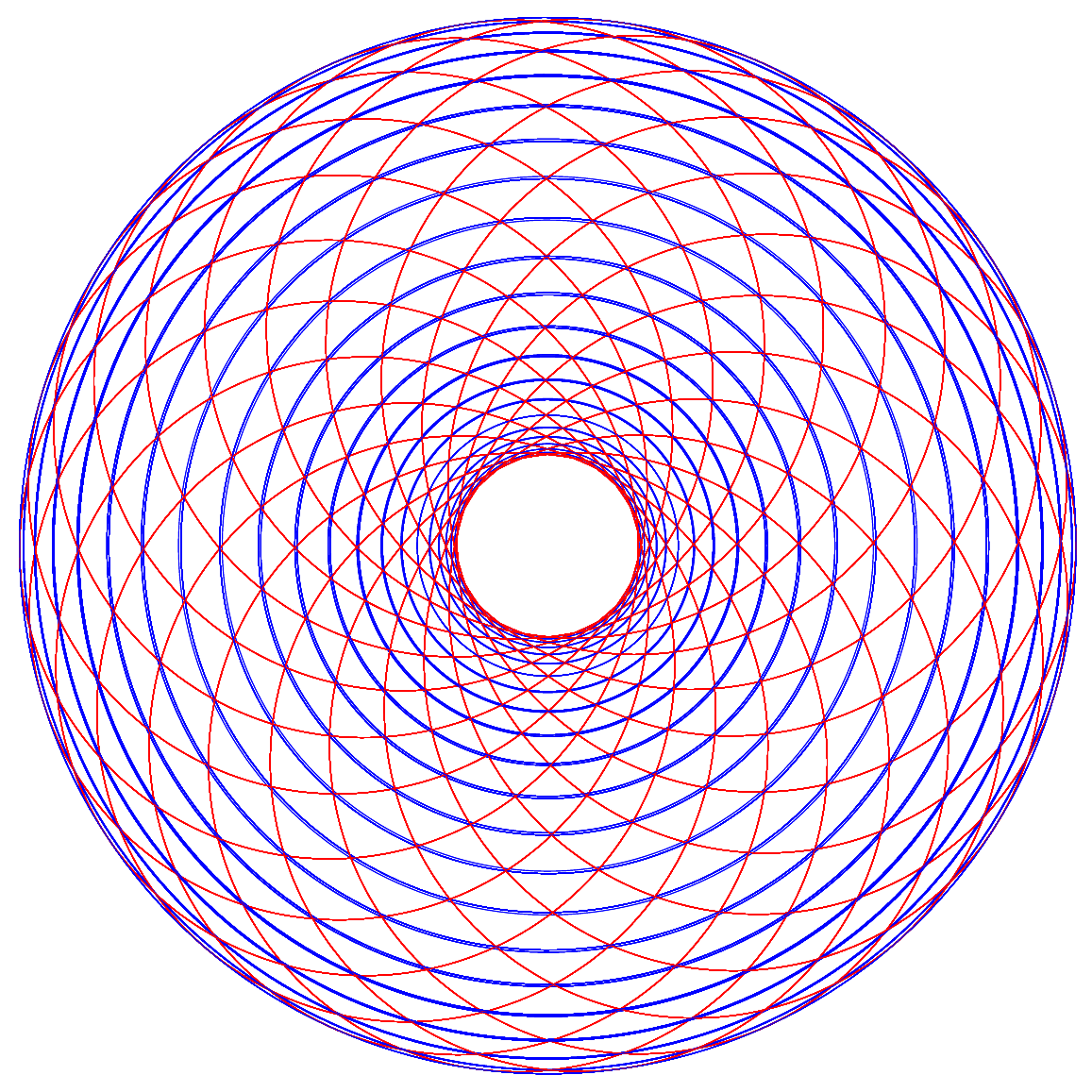}
\includegraphics[width=0.30\textwidth, trim={0.3cm 3.5cm 8cm 3.5cm}, clip]{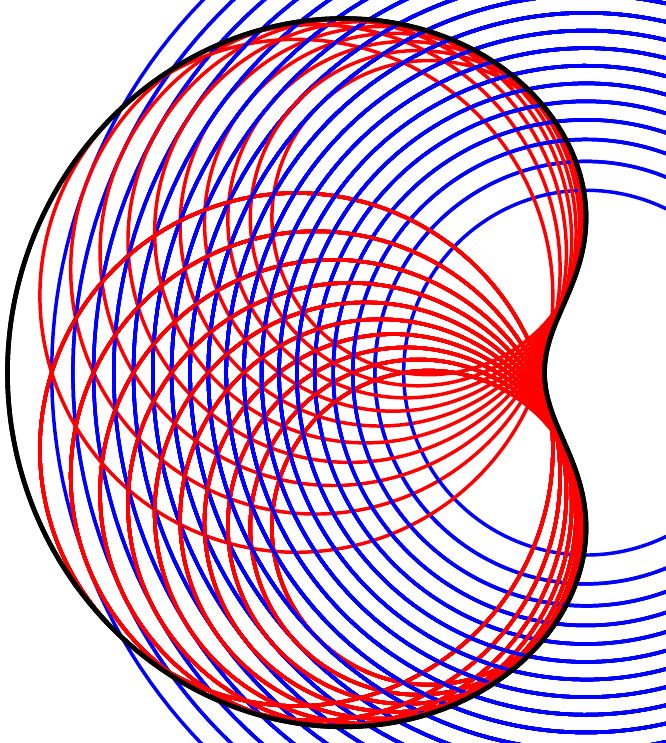}
 \caption{Types 1, 2, 3.}\label{123}
\end{figure}

 {\bf 1. Polar conic plane does not cut Darboux quadric, hyperbolic pencil.}
 The pencil with limit circles at the origin and in the infinite point gives circles $x^2+y^2=v$, the polar conic is the circle
$
 X_0^2 + Y_0^2+4cX_0Z_0 =0$, $ U_0=0$, $ c> 0$,
 defining the family
 \[
 x^2+ y^2 -\frac{4c}{c^2u^2 + 1}x +\frac{4c^2u}{c^2u^2 + 1}y=1,
 \]
 the circles of the family enveloping the cyclic
 \[ \bigl(x^2 + y^2\bigr)^2 - \bigl(x^2 + y^2\bigr)(4cx + 2) -4c^2y^2 + 4cx + 1=0
 \]
 as shown on the left in Figure~\ref{123}.

 {\bf 2. Polar conic plane does not cut Darboux quadric, hyperbolic pencil, webs symmetric by rotations.}
 The pencil with limit circles at the origin and at the infinite point gives circles $x^2+y^2=v$, the polar conic is the circle
\smash{$
 X_0^2 + Y_0^2 = \frac{4Z_0^2}{c^2}$}, $ U_0=0$,
 defining the family
 \[
 x^2+ y^2 + \frac{2\cos(u)}{c}x+\frac{2\sin(u)}{c} y=1,
 \]
 the circles of the family enveloping the cyclic{\samepage \[
 c^2\bigl(x^2 + y^2\bigr)^2 - \bigl(2c^2 + 4\bigr)\bigl(x^2 + y^2\bigr) + c^2=0,
 \]
 which splits into two concentric circles as shown in the center of Figure~\ref{123}.}

{\bf 3. Polar conic plane cuts Darboux quadric, hyperbolic pencil.}
 The pencil with limit circles at the origin and at the infinite point gives circles
$
x^2+y^2=v$,
the polar conic is the circle
$
x_0^2 + y_0^2=4cx_0$, $ z_0=0$, $ c>0$,
 defining the family
 \[
 x^2+ y^2 +\frac{4c}{c^2u^2 + 1} x - \frac{4c^2u}{c^2u^2 + 1}y=-1,
 \]
 the circles of the family enveloping the cyclic \[
 \bigl(x^2 + y^2\bigr)^2 + \bigl(x^2 + y^2\bigr)(4cx + 2) - 4c^2y^2 + 4cx + 1=0,
 \]
 as shown on the right of Figure~\ref{123}.

\begin{figure}[t]\centering
\includegraphics[width=0.3\textwidth, trim={0.4cm 0.3cm 0.3cm 0.3cm}, clip]{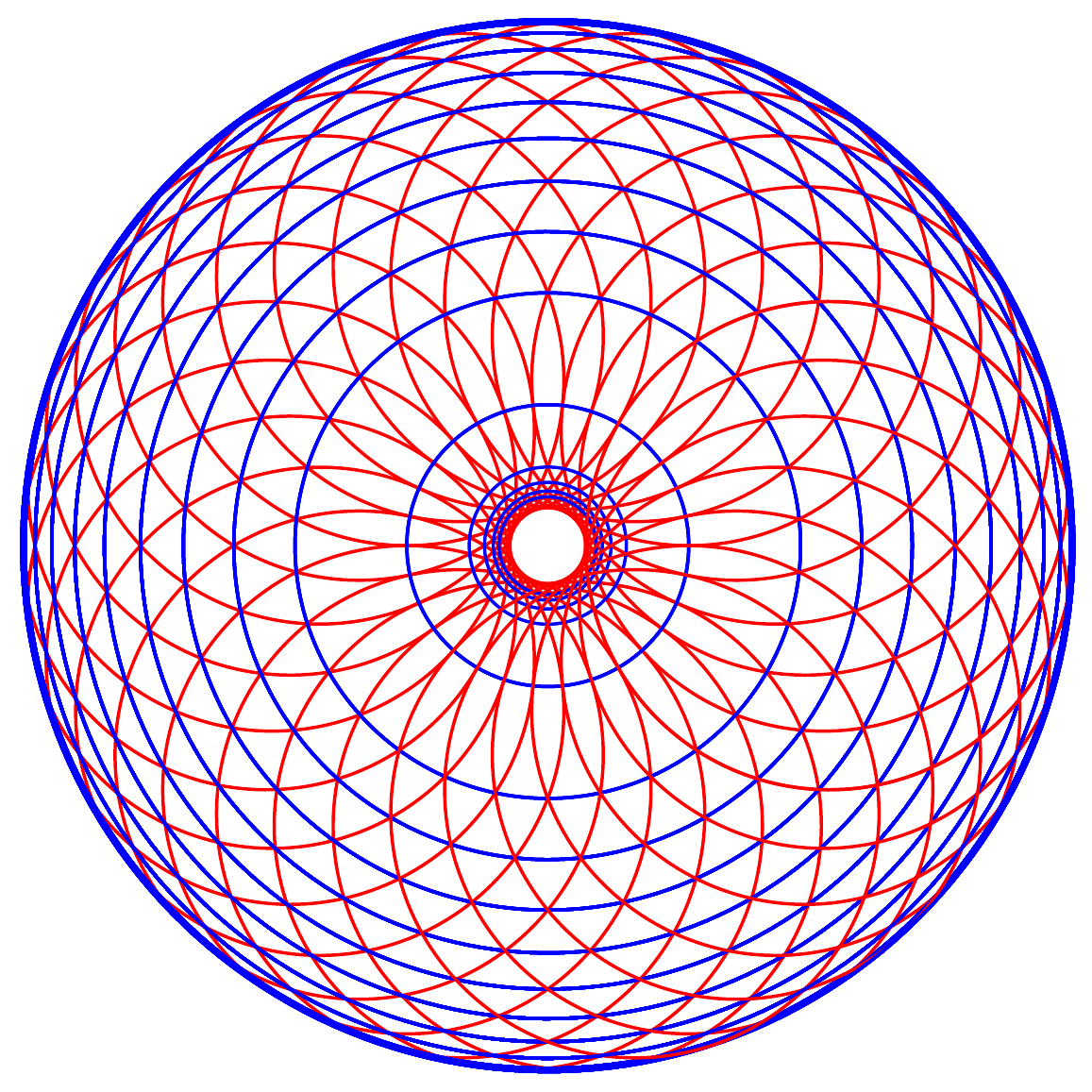}
\includegraphics[angle=0, width=0.33\textwidth, trim={0.5cm 3.0cm 0.0cm 0.5cm}, clip]{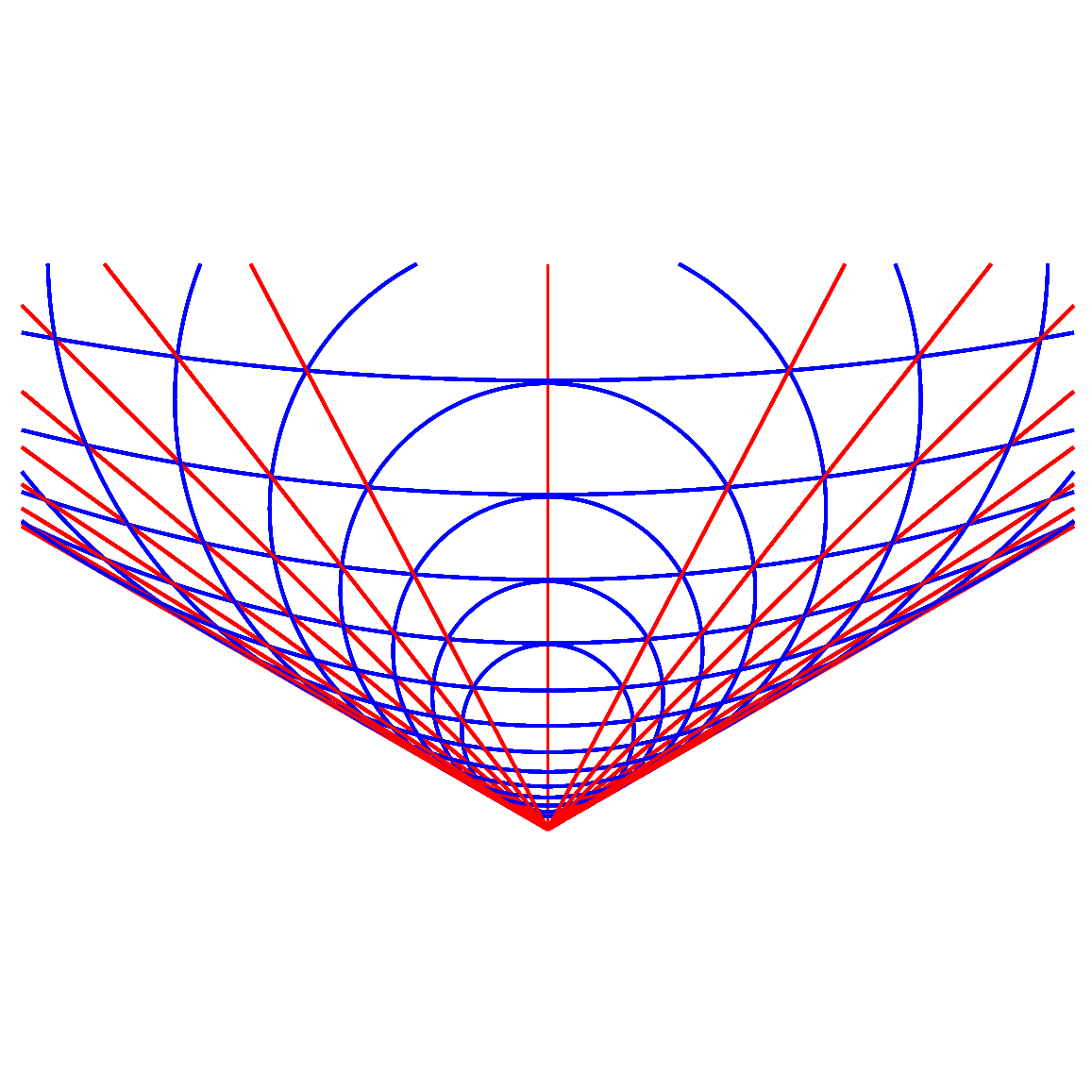}
\includegraphics[angle=0, width=0.3\textwidth, trim={0.4cm 0.0cm 0.0cm 0.0cm}, clip]{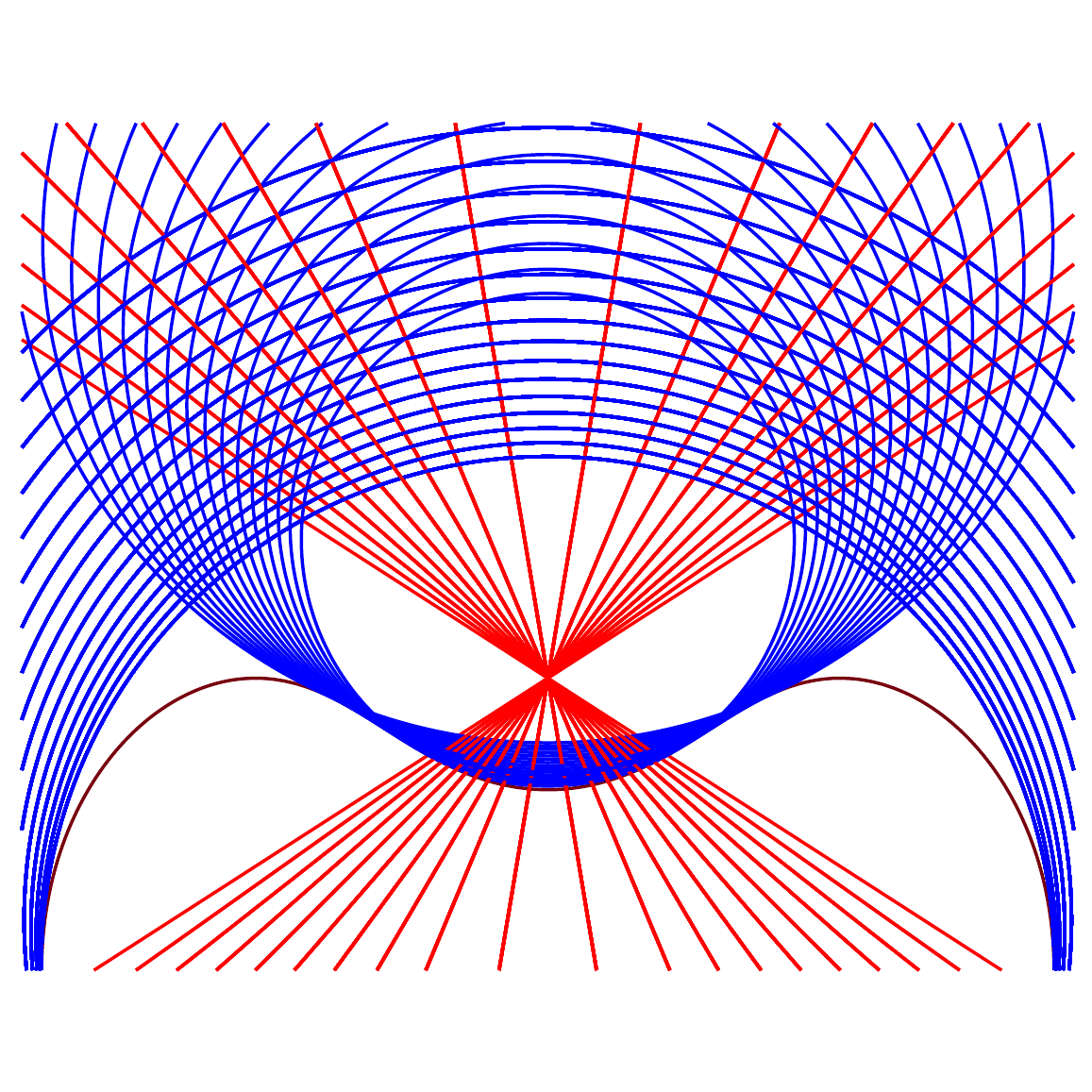}
 \caption{Types 4, 5, 6.}\label{456}
\end{figure}

 {\bf 4. Polar conic plane cuts Darboux quadric, hyperbolic pencil, webs symmetric by rotations.}
 The pencil with limit circles at the origin and at the infinite point gives circles~${
x^2+y^2=v}$,
the polar conic is the circle
$
x_0^2 + y_0^2=\frac{4}{c^2}$, $ z_0=0$,
 defining the family
 \[
 x^2+ y^2 -\frac{2\cos(u)}{c}x-\frac{2\sin(u)}{c} y=-1,
 \]
 the circles of the family enveloping the cyclic \[
 c^2\bigl(x^2 + y^2\bigr)^2 + \bigl(2c^2 - 4\bigr)\bigl(x^2 + y^2\bigr) + c^2=0,\]
 which splits into two concentric circles as shown on the left in Figure~\ref{456}.

 {\bf 5. Polar conic plane cuts Darboux quadric, elliptic pencil, webs symmetric by homotheties.}
 The pencil with vertexes at the origin and at the infinite point gives lines
$
y=vx$,
the polar conic is
$
y_0^2+ cz_0^2= c$, $ x_0=0$, $ c\notin [0,1]$,
defining the family of circles
$
x^2 + (y + u)^2 =(1-1/c)u^2$,
the circles enveloping the lines $x^2=(c-1)y^2$,
real for $c>1$, as shown in the center in Figure~\ref{456}. The webs of the family are symmetric by homotheties with the center in the origin.

 {\bf 6. Polar conic plane cuts Darboux quadric, elliptic pencil, polar conic and dual to polar line intersect in 2 points on the Darboux quadric.}
 The pencil with vertexes at the origin and at the infinite point gives lines $y=vx$, the polar conic is
$
2cy_0=z_0^2-1$, $ x_0=0$, $ c>0$,
 defining the family of circles
\[
x^2 + y^2 +\frac{(u - 1)}{cu} y + \frac{u - 1}{u + 1}=0,
\]
the circles enveloping the cyclic
\[
c^2\bigl(x^2 + y^2\bigr)^2 + \bigl( 4cy-2c^2 \bigr)\bigl(x^2 + y^2\bigr) + (2y+c)^2=0\]
 as shown on the right in Figure~\ref{456}.

 \begin{figure}[t]\centering
\includegraphics[width=1\textwidth]{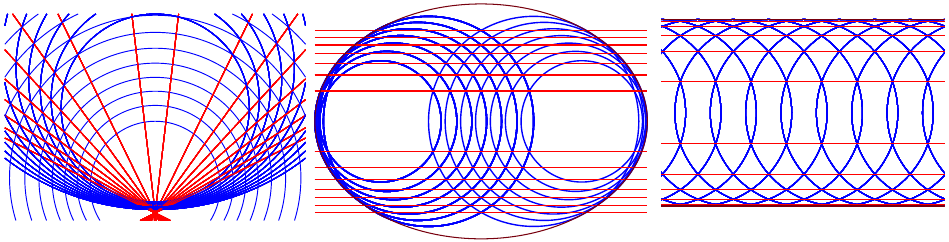}
 \caption{Types 7, 8, 9.}\label{789}
\end{figure}

 {\bf 7. Polar conic plane cuts Darboux quadric, elliptic pencil.}
 The pencil with vertexes at the origin and at the infinite point gives lines
$y=vx$,
 the polar conic is
$
 2cy_0z_0+1=z_0^2$, $ x_0=0$, $ c>0$,
 defining the family of circles
\[
x^2 + y^2 - \frac{c+u}{cu}y + \frac{c + u}{c - u}=0,
\]
 as shown on the left in Figure~\ref{789}.

 {\bf 8. Polar conic plane cuts Darboux quadric, parabolic pencil.}
 The pencil with the vertex at the infinite point gives lines $y=v$, the polar conic
$
x_0^2 - z_0 =1$, $ y_0=0$,
 gives the family of circles
\smash{$
\big(x - \frac{1}{u}\big)^2 + y^2 = \frac{u^2 - 1}{u^2}$}.
The circles envelopes the ellipse $x^2 + 2y^2 = 2$ as shown in the center in Figure~\ref{789}.

 {\bf 9. Polar conic plane cuts Darboux quadric, parabolic pencil, webs symmetric by translations.}
 The pencil with the vertex at the infinite point gives lines $y=v$, the polar conic~is
$
x_0^2 +\bigl(1-\sqrt{2}\bigr)z_0^2 - 2\sqrt{2}z_0 = 1+\sqrt{2}$, $ y_0=0$,
defining the family of circles
$
(x +u)^2 + y^2 =1$,
the circles touching the lines $y=\pm 1$ as shown on the right in Figure~\ref{789}.

\begin{figure}[t]\centering
\includegraphics[width=0.30\textwidth, trim={9.2cm 6.35cm 2.8cm 6.4cm}, clip]{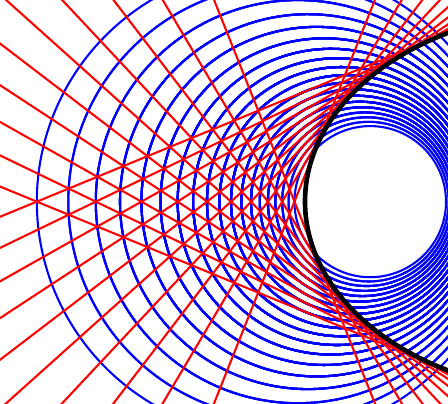}
\includegraphics[width=0.30\textwidth, trim={3.1cm 3.5cm 3.5cm 3.5cm}, clip]{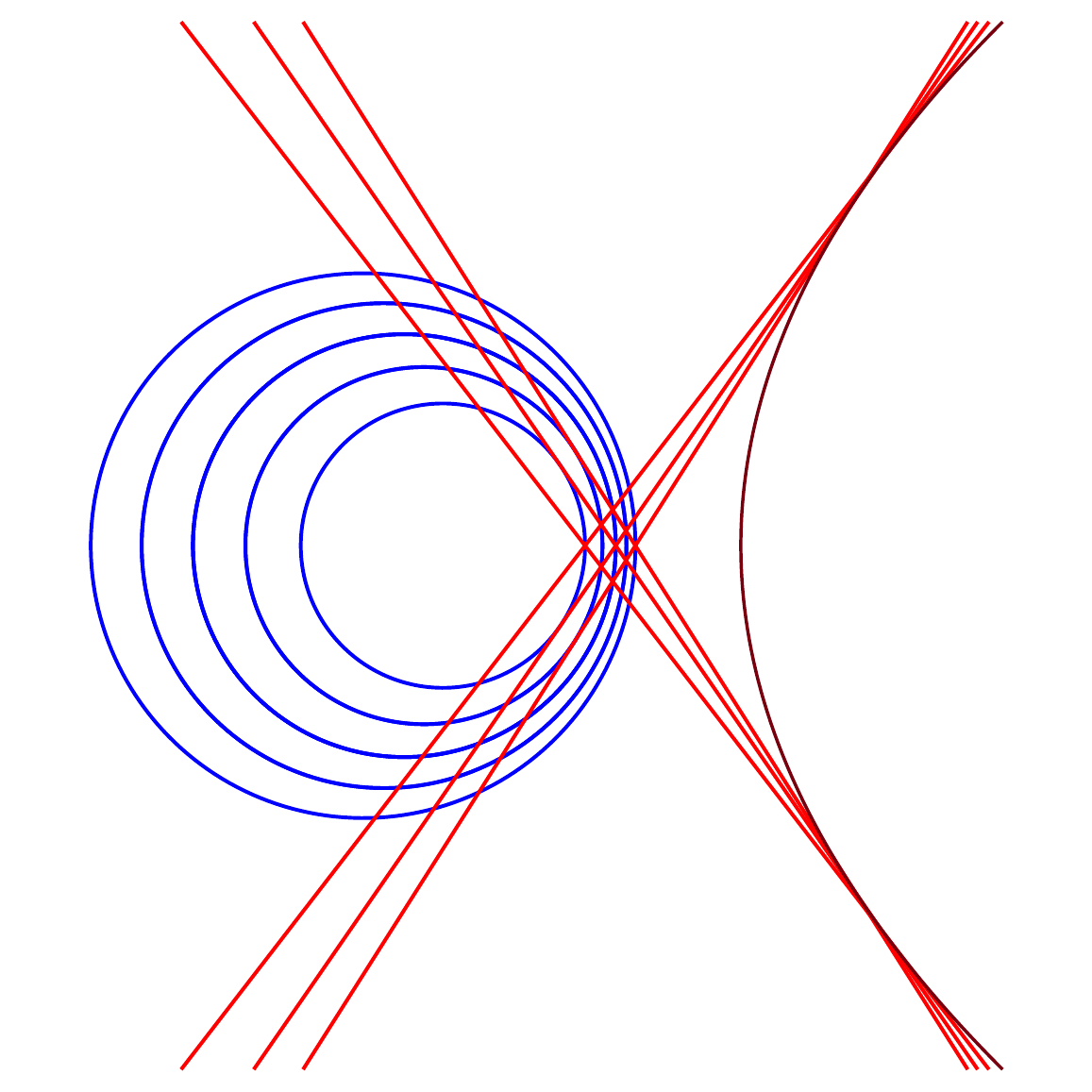}
\includegraphics[width=0.36\textwidth, trim={3.3cm 6cm 3.5cm 4cm}, clip]{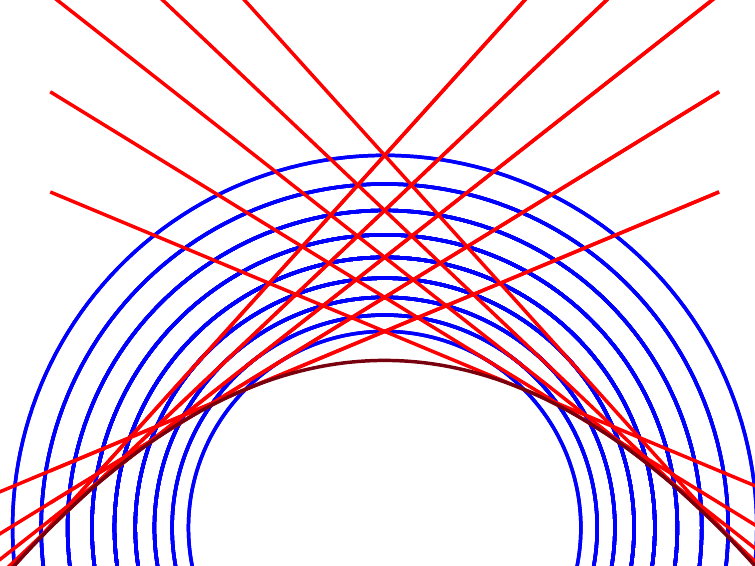}
 \caption{Type 10, 11, 12.}\label{101112}
\end{figure}

 {\bf 10. Polar conic plane tangent to Darboux quadric, hyperbolic pencil.}
 The pencil with limit circles at $(1,0)$ and $(-1,0)$ gives circles
$
x^2 + y^2+ 1=vx$,
 the polar conic is
$
x_0^2+ {(1-c)y_0^2} =c$, $ z_0=-1$, $ c>0$, $ c\ne 1$,
 defining the family of lines $x_0x + y_0y = 1$, the lines of the family enveloping the conic
$
 cx^2 + \frac{c}{c - 1}y^2=1
$
 with foci $(1,0)$ and $(-1,0)$ as shown on the left in Figure~\ref{101112}.

 {\bf 11. Polar conic plane tangent to Darboux quadric, hyperbolic pencil, polar line meets polar conic.}
 The pencil with limit circles at $(1,0)$ and $(-1,0)$ gives circles
$
x^2 + y^2+ 1=vx$,
the polar conic is
$
y_0^2 -cx_0y_0 + cy_0 + x_0=0$, $ z_0=-1$, $ c\ge 0$,
 defining the family of lines~${
x_0x + y_0y = 1}$,
the lines of the family enveloping the parabola{\samepage
\[
(cx - y)^2 - (2c^2 + 4)x - 2cy + c^2=0
\] with focus at $(1,0)$ as shown in the center in Figure~\ref{101112}.}

 {\bf 12. Polar conic plane tangent to Darboux quadric, hyperbolic pencil, polar line contains the point where polar conic plane touches the Darboux quadric.}
 The pencil with limit circles at the origin and infinity gives circles
$
x^2 + y^2=v$,
 the polar conic is
$
x_0^2 + y_0^2 - 2y_0=0$, $ z_0=-1$
defining the family of lines
$
x_0x + y_0y = 1$,
the lines of the family enveloping the parabola $ x^2 + 2y = 1$ with focus at $(0,0)$ as shown on the right in Figure~\ref{101112}.

\begin{figure}[t]\centering
\includegraphics[width=0.315\textwidth, trim={0.5cm 0.5cm 0.5cm 0.5cm}, clip]{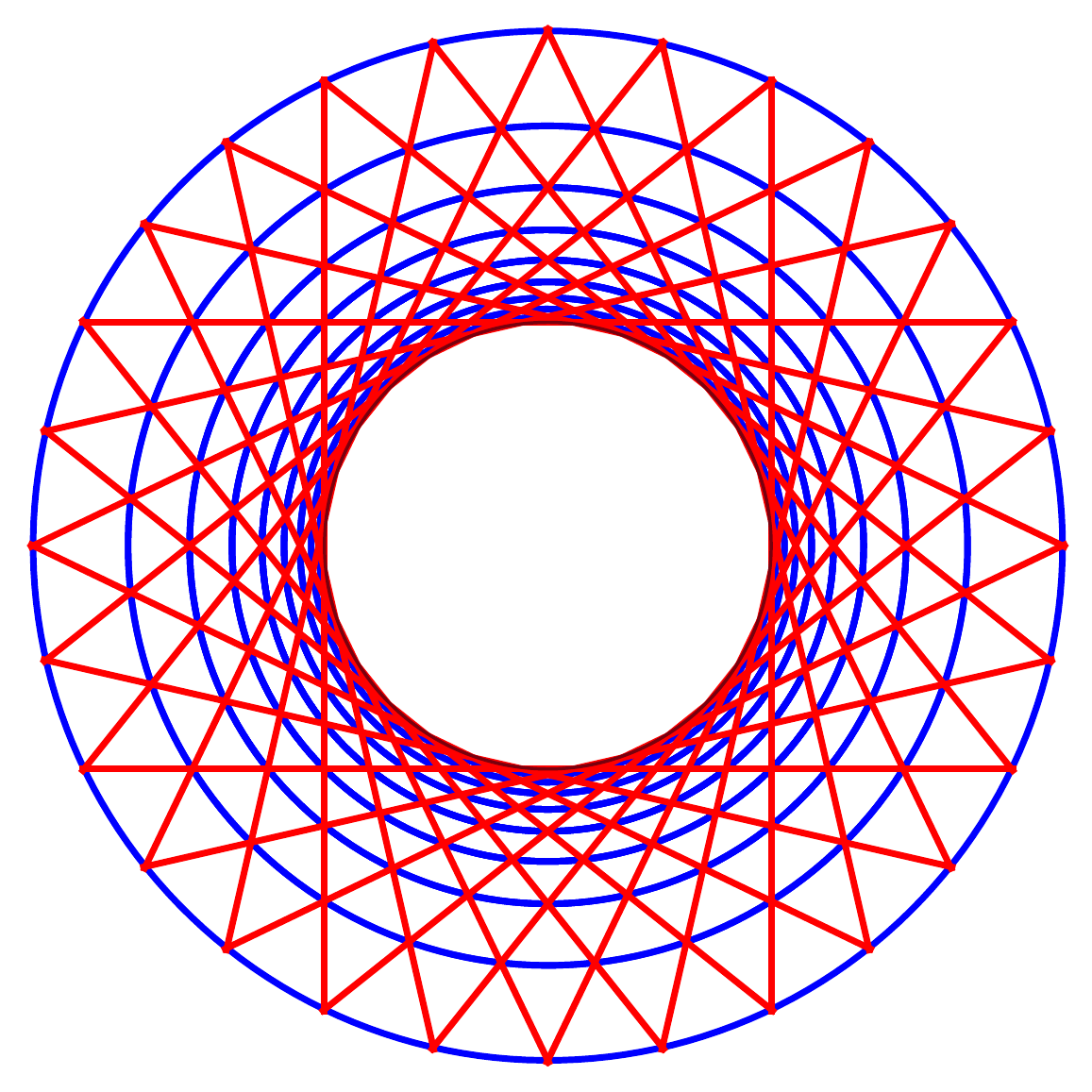}
\includegraphics[width=0.335\textwidth, trim={5cm 3.2cm 5cm 7.8cm}, clip]{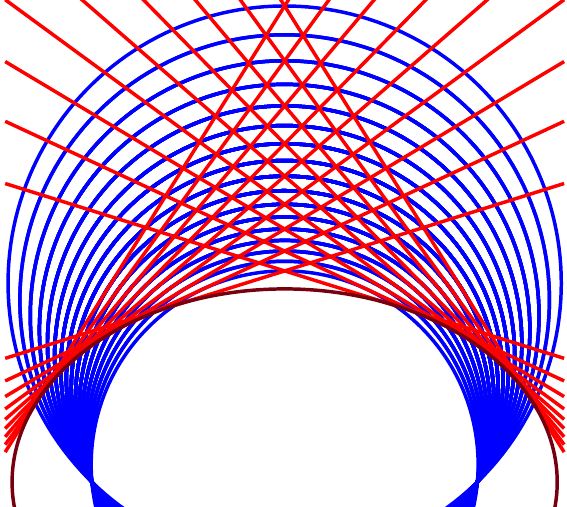}
\includegraphics[width=0.325\textwidth, trim={5.5cm 5cm 5.5cm 6.5cm}, clip]{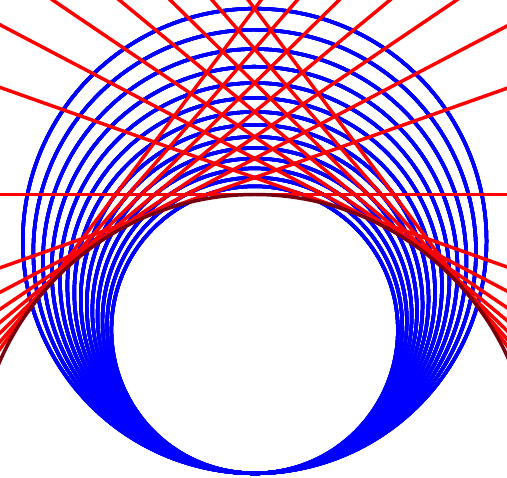}
 \caption{Type 13, 14, 15.} \label{131415}
\end{figure}

 {\bf 13. Polar conic plane tangent to Darboux quadric, hyperbolic pencil, web symmetric by rotations.}
 The pencil with limit circles at the origin and infinity gives circles $x^2 + y^2=v$, the polar conic is $x_0^2 + y_0^2 =1$, $ z_0=-1$ defining the family of lines $x_0x + y_0y = 1$, the lines of the family enveloping the circle $ x^2 + y^2 = 1$ as shown on the left in Figure~\ref{131415}.

 {\bf 14. Polar conic plane tangent to Darboux quadric, elliptic pencil.}
 The pencil with vertexes at $(1,0)$ and $(-1,0)$ gives circles
\[
\frac{\bigl(x^2 + y^2-1\bigr)^2}{\bigl(x^2 + y^2-2x+1\bigr)\bigl(x^2 + y^2+2x+1\bigr)}=v,
\]
the polar conic is
$
cx_0^2+ (c + 1)y_0^2 + 1=0$, $ z_0=-1$,
 defining the family of lines
$
x_0x + y_0y = 1$,
 the lines of the family enveloping the conic
\smash{$
 \frac{x^2}{c} + \frac{y^2}{c + 1}=-1
 $} with foci $(1,0)$ and $(-1,0)$ as shown in the center in Figure~\ref{131415}.

 {\bf 15. Polar conic plane tangent to Darboux quadric, parabolic pencil.}
 The pencil with vertex at the origin gives circles
$
x^2 + y^2=2vy$,
the polar conic is
$
x_0^2 + y_0^2 =1$, $ z_0=-1
$ defining the family of lines $x_0x + y_0y = 1$, the lines of the family enveloping the circle $ x^2 + y^2 = 1$ as shown on the right in Figure~\ref{131415}.

\begin{Theorem}
The webs of different types in the above classification list are not M\"obius equivalent. The webs of the same type with different normal forms are not M\"obius equivalent.
\end{Theorem}
\begin{proof} The webs from different types are not M\"obius equivalent: discrete geometric invariants indicated in the descriptions, such as 1) presence of infinitesimal symmetry and 2) the mutual position of polar line, polar conic and Darboux quadric, effectively separate the types.

To see that the different normal forms within a family are not M\"obius equivalent, one computes the subgroup $G_p$ of ${\rm PSO}(3,1)$ respecting the positions of the polar conic plane and the polar line in the chosen normalization and checks that the $G_p$-orbits of the canonical forms are different.

For the first 4 types, $G_p$ is generated by $R_z$ and by reflections in the coordinate planes.

For the 5th, 6th, 7th type, $G_p$ is generated by $B_z$ and by reflections in the coordinate planes.

For the 10th, 11th and 14th types, $G_p$ is discrete and generated by reflections in the planes $x=0$ and $y=0$.
\end{proof}

As in the case of 3 pencils, Lemma \ref{sing} effectively selects candidates among 3-webs that can be hexagonal. The singularities of the type described by the Lemma arise when either 1) a line, joining 2 different points on the polar conic, touches the Darboux quadric at a point $p_t$, while the polar conic plane is not tangent the Darboux quadric or 2) a line, joining a point on the polar line with a point on the polar conic, touches the Darboux quadric at a point $p_t$ while the tangent plane to the Darboux quadric at $p_t$ is not tangent to the polar conic.

In the former case, the points $p_t$ trace a circle, which is the intersection of the polar conic plane with the Darboux quadric. Then, by Lemma \ref{sing}, the polar of the polar conic plane lies on the pencil polar line.

In the latter case, consider a point $p_r$ running over the polar line. For non-planar polar set, $p_r$ meets the polar conic plane $\pi_c$ only at one point. Therefore, the one-parameter family of cones tangent to the Darboux quadric and having their vertexes at $p_r$ cuts the plane $\pi_c$ in a one-parameter family of conics $c_r$. Each conic $c_r$ intersect the polar conic $c_p$ at 4 (possibly complex or multiple) points. If the polar conic $c_p$ is not a member of the family $\{c_r\}$, at least one of these 4 intersection points is moving along $c_p$ as $p_r$ runs over the polar line. In fact, if intersection points are stable then $\{c_r\}$ is a pencil of conics containing $c_p$.

Choose one such moving intersection point $p_i$. Lines $l_r$ tangent to $c_r$ at $p_i$ form a one-parameter family, or congruence of lines. All the objects in this construction are considered as complex but the polar conic and the polar line must have equations with real coefficients and the real part of the polar conic cannot lie completely inside the Darboux quadric.

\begin{Proposition}\label{geometryconic}
If the polar curve of a hexagonal circular $3$-web is non-planar and splits into a line and a smooth conic, then for the web complexification hold true
\begin{itemize}\itemsep=0pt
\item[$(1)$] the polar of the polar conic plane lies on the polar line, if this plane is not tangent to the Darboux quadric and
\item[$(2)$] the congruence of lines $l_r$ is a pencil with the vertex on the polar conic, if the polar conic is not a member of the family $\{c_r\}$.
\end{itemize}
\end{Proposition}
\begin{proof}
 The first claim follows directly from Lemma \ref{sing}: the points, where bisecant lines touch the Darboux quadric, trace a circle whose polar point must lie on the polar line. For conic planes missing the Darboux quadric, the complex version works.

To derive the second claim, observe that the line $p_rp_i$ touches the Darboux quadric at a~singular point treated by Lemma \ref{sing}. Thus the tangency point must trace a circle on the Darboux quadric and the polar point $p_0$ of the circle must lie on the polar conic. Then the plane, tangent to the Darboux quadric and passing through $p_rp_i$, contains $p_0$. This plane cuts the plane $\pi_c$ of the polar conic along a line $l_r$, passing through the $p_i$ and tangent to the corresponding conic~$c_r$. Hence $\{l_r\}$ is the pencil with vertex at $p_0$.
\end{proof}

\begin{Theorem}\label{conicPencil}
If the polar curve of a hexagonal circular $3$-web splits into a smooth conic and a straight line, not lying in a plane of the conic, then the web is M\"obius equivalent to one from the above presented list of types $1$--$15$.
\end{Theorem}
\begin{proof} The polar conic plane either completely misses the Darboux quadric, or cuts it in a~circle, or is tangent to it. The polar line is either hyperbolic, or elliptic, or parabolic. Thus we have 9 cases to consider, each case defining a set of webs (possibly empty).

$\bullet$ {\it Polar conic plane misses the Darboux quadric}.
Applying a suitable M\" obius transformation, we can send the polar conic plane to infinity. Then by Proposition \ref{geometryconic} the polar line contains the origin of the affine chart and therefore meets the Darboux quadric at 2 points. Thus the polar line is hyperbolic. Applying a rotation around the origin, we map these two points to~${(0,0,\pm 1)}$. In the affine coordinates $x=\frac{X}{Z}$, $ y=\frac{Y}{Z}$ on the polar conic plane $U=0$, the conics~$c_r$ are~${
x^2+y^2=r}$,
where $r$ is considered as a complex parameter. These conics are real only for real non-negative $r$.
If the polar conic coincides with one of the conics~$c_r$, then we get Type~2, symmetric by~$R_z$.

If the polar conic is not one of $c_r$, then, by Proposition \ref{geometryconic}, the points, where lines of pencil with the vertex at some point $p_0=(x_0,y_0)\in \mathbb C^2$ touch the circles $c_r$, run over the polar conic. The line of the pencil $y=y_0+k(x-x_0)$ corresponding to the parameter $k\in \mathbb C$, is tangent to the conic $c_r$ for
\[
r=\frac{x_0^2k^2 - 2x_0y_0k + y_0^2}{k^2 + 1}
\]
at the point
\[
p(k)=(x(k),y(k))=\left(\frac{k(x_0k - y_0)}{k^2 + 1},\frac{y_0-x_0k}{k^2 + 1} \right).
\]
The points $p(k)$ run over the complex conic
$
x^2+y^2 =x_0x + y_0x$.
 This conic is real if and only if~$x_0$ and $y_0$ are real. It is smooth if and only if $p_0\ne (0,0)$. We got the first web of the list, Type~1. This conic is the circle, passing through the origin $O=(0,0)$ and $p_0$ and having its center at the midpoint of the segment $Op_0$.
(We obtained a theorem of scholar geometry.) Using~$R_z$, we normalize $p_0$ to $y_0=0$, $x_0>0$.

 For infinite $p_0$, different from the cyclic points, the line $l(k)$ of the pencil $y=ax+k$, $a\in \mathbb C$ touches a unique conic of the family $\{c_r\}$ at the point
$
p(k)=\bigl(-\frac{ak}{a^2 + 1}, \frac{k}{a^2 + 1}\bigr)$.
The points $p(k)$ are collinear: $x(k)+ay(k)=0$ and we cannot obtain a smooth polar conic in this way. Finally, if $p_0$ is cyclic, the lines from the pencil touch the conics $c_r$ at the very point~$p_0$ and we get no conic at all.

$\bullet$ {\it Polar conic plane cuts the Darboux quadric, hyperbolic pencil}.
We send the conic plane to~${Z=0}$. Now the polar line contains the infinite point $[0:0:1:0]$ by Proposition \ref{geometryconic}. The subgroup of the M\" obius group, preserving the plane $Z=0$, is generated by rotations around $Z$-axis and the boosts $B_x$ and $B_y$. Using this subgroup, we normalize the polar line to~${x=y=0}$. The conics $c_r$ are circles in the plane $z=0$ with the center at the origin. Repeating the arguments that we used above for conic planes missing the Darboux quadric, we get Types 3 and 4.

$\bullet$ {\it Polar conic plane cuts the Darboux quadric, elliptic pencil}.
We normalize the conic plane to $X=0$. Then by Proposition \ref{geometryconic} the point $[1:0:0:0]$ lies on the polar line. The polar line, being elliptic, meets the plane $X=0$ at a point $p$ outside the unit circle. M\" obius transformations, preserving the plane $X=0$, are generated by rotations around $X$-axis and the boosts $B_y$ and $B_z$. Using these transformations, we send the point $p$ to $[0:1:0:0:0]$. Now the polar line is $U=Z=0$ and the conics $c_r$ have equations
\smash{$
z^2+\frac{ y^2}{r}=1
$} in the affine coordinates.
If the polar conic coincides with one of these conics, then $r$ is real and we get Type 5.

 Otherwise, by the second claim of Proposition \ref{geometryconic}, the lines $l_r$ meet at some point $p_0\in c_p$. If $p_0=(0,y_0,z_0)$ is finite, a line of the pencil of lines $z=z_0+k(y - y_0)$ with vertex at $p_0$ is tangent to the conic $c_r$ with
\[
r=\frac{y_0^2k^2 - 2y_0z_0k + z_0^2 - 1}{k^2}.
\]
 Thus the points, where the lines $l_r$ touch $c_r$, are parametrized by $k$ via
\[
y=\frac{y_0^2k^2 - 2y_0z_0k + z_0^2 - 1}{k(y_0k - z_0)}, \qquad z=\frac{1}{z_0-y_0k}.
\]
Excluding $k$, we get the conic
$
y_0z^2-z_0yz + y - y_0=0$.
This conic is real only for real $y_0$, $z_0$ and is smooth if and only if $y_0\bigl(z^2_0-1\bigr)\ne 0$.
 If $z^2_0<1$, we normalize to $z_0=0$ applying $B_z$ and get Type~6. If $z^2_0> 1$, we send $p_0$ to infinity applying~$B_z$.

For infinite point $p_0=[0:Y_0:Z_0:0]$, we conclude that $Y_0\ne 0$ and $Z_0\ne 0$, otherwise the tangency points $p_i$ do not trace a conic. Thus we can set $p_0=[0:y_0:1:0]$, where $y_0\ne 0$. In~the affine coordinates $u=\frac{U}{Z}$, $y=\frac{Y}{Z}$ in the plane $X=0$, the conic $c_r$ is \smash{$1+\frac{y^2}{r}=u^2$}. A~line from the pencil $y=ku+y_0$ is tangent to the conic $c_r$ if and only if
$
r=k^2 - y_0^2$.
The points, where the lines $l_r$ touch $c_r$, are parameterized by $k$ via
\smash{$
y=\frac{y_0^2 - k^2}{y_0}$}, \smash{$ u=-\frac{k}{y_0}$}.
This is a parametrization of the polar conic of Type 7
$
y_0u^2 + y = y_0$,
which becomes
$
y_0U^2 - y_0Z^2 + YZ=0
$
in homogeneous coordinates.

$\bullet$ {\it Polar conic plane cuts the Darboux quadric, parabolic pencil}.
We normalize the conic plane to $Y=0$. Then the point $[0:1:0:0]$ is on the polar line.
 The plane $Y=0$ is stable under rotations $R_y$ along the $y$-axis and under boosts $B_x$ and $B_z$.
Rotating, if necessary, around the $y$-axis, we bring the polar line to $z=-1,$ $x=0.$ In the affine coordinates on the plane $Y=0$, the conics $c_r$ have equations
\[
x^2 + \frac{r-1}{r}z^2 - \frac{2}{r}z - \frac{r+1}{r}=0.
\]
If the polar conic coincides with one of $c_r$, then $r$ is real. The chosen position of polar conic plane and polar line is stable by action of the group generated by $B_z$. Applying it, one can normalize to $r=\frac{1}{\sqrt{2}}$ and we get Type 9.

Otherwise, consider first the pencil of lines $z=z_0+k(x - x_0)$ with finite vertex at $p_0=(x_0,0,z_0)$. A line from the pencil is tangent to the conic $c_r$ with
\[
r=\frac{x_0^2k^2 - 2x_0(z_0 + 1)k + (z0 + 1)^2}{\bigl(x_0^2 - 1\bigr)k^2 - 2x_0z_0k + z_0^2 - 1}.
\]
 Thus the points, where the lines $l_r$ touch $c_r$, are parametrized by $k$ via
\[
x=\frac{k(x_0k - z_0 - 1)}{k^2 - x_0k + z_0 + 1}, \qquad z=\frac{\bigl(x_0^2 - 1\bigr)k^2 - x_0(2z_0 + 1)k + z_0(z_0 + 1)}{k^2 - x_0k + z_0 + 1}.
\]
Excluding $k$, we get the conic
\begin{equation}\label{polarcutP}
(z_0 + 1)x^2 - x_0xz +z^2- x_0x+ (1-z_0)z - z_0=0.
\end{equation}

This conic is real only for real $x_0,z_0$ and is smooth if and only if $z_0\ne -1.$
The chosen position of polar conic plane and polar line is stable by action of the group generated by $B_z$ and $B_x-R_y$. Consider the orbits of points in the plane $Y=0.$ The orbit dimension is two for points outside the union of the line $U+Z=0$ and the circle $X^2+Z^2=U^2$, and is one on this union except for their common point.
Thus the finite representatives of the orbits are $[0:0:0:1]$, $[0:0:1:1]$ and $[0:0:-1:1].$ For the first two points, the conics \eqref{polarcutP} lie inside the Darboux quadric and there is no real circles. For the point $[0:0:-1:1]$, the conic \eqref{polarcutP} is not smooth.

For infinite vertexes $p_0$, the pencil $z=k$ gives a non-smooth conic. Thus the pencil can be chosen as $x=az+k$. A line from the pencil is tangent to the conic $c_r$ with
\smash{$
r=\frac{(k-a)^2}{k^2 -a^2 - 1}$}.
 The points, where the lines $l_r$ touch $c_r$, are
\smash{$
x=\frac{k-a}{a^2 - ak + 1}$}, \smash{$ z=\frac{k^2-ak-1}{a^2 - ak + 1}$}.
Excluding $k$, we get the conic~${
x^2-axz-ax-z-1=0}$,
which is real only for real $a$. Taking into account the action of the group generated by $B_z$ and $B_x-R_y$, we set $a=0$ and get Type 8.

$\bullet$ {\it Polar conic plane tangent to Darboux quadric, hyperbolic pencil}.
 If the hyperbolic line does not contain the point where the polar conic plane touches the Darboux quadric, then we sent this point to $(0,0,-1)$ and the polar line to $Y=Z=0$. In the affine coordinates on the plane~${z=-1}$, the conics $c_r$ have equations
$(x - r)^2 + \bigl(1-r^2\bigr)y^2 + 1-r^2=0$.
If the polar conic coincides with one of $c_r$, then the web curvature
\[
K_B=\frac{4\bigl(r^2 - 1\bigr)^2\bigl(x^2 - 1\bigr)\bigl( rx^2 + ry^2 +\bigl(r^2-3\bigr)x +r\bigr)\bigl(x^2 + y^2-2rx + 1\bigr)^4}{x^4y^3\bigl(x^2 +1-2rx \bigr)^6}
\]
vanishes only for $r=\pm 1$, the conic $c_r$ being non-smooth for these values.

A line from the pencil $y=y_0+k(x - x_0)$ with finite vertex at $p_0=(x_0,y_0,-1)$ is tangent to the conic $c_r$ with
\[
r=\frac{\bigl(x_0^2 + 1\bigr)k^2 - 2x_0y_0k + y_0^2 + 1}{2k(x_0k - y_0)}.
\]
 Thus the points, where the lines $l_r$ touch $c_r$, are parametrized by $k$ via
\begin{gather*}
x=\frac{x_0\bigl(x_0^2 - 1\bigr)k^3 - y_0\bigl(3x_0^2 - 1\bigr)k^2 + x_0\bigl(3y_0^2 + 1\bigr)k - y_0\bigl(y_0^2 + 1\bigr)}{k\bigl(\bigl(x_0^2 - 1\bigr)k^2 - 2x_0y_0k + \bigl(y_0^2 - 1\bigr)\bigr)},
\\
y=\frac{2(x_0k - y_0)}{\bigl(x_0^2 - 1\bigr)k^2 - 2x_0y_0k + \bigl(y_0^2 - 1\bigr)}.
\end{gather*}
Excluding $k$, we get the cubic
\begin{equation}\label{polartH}
\bigl(x_0^2 - 1\bigr)y^3 +\bigl(y_0^2 - 1\bigr)x^2y - 2x_0y_0xy^2+ 2y_0x^2+ 2y_0y^2- 2x_0y_0x +\bigl(x_0^2 - y_0^2\bigr)y=0.
\end{equation}
This cubic splits into a smooth real conic and a line in 3 cases:
\begin{itemize}\itemsep=0pt
\item[(1)] for $y_0=0$, \eqref{polartH} factors as $y\bigl(x^2 +\bigl(1-x_0^2\bigr)y^2 - x_0^2\bigr)=0$ and we get Type 10,
\item[(2)] for $x_0= 1$, $y_0\ne 0$, \eqref{polartH} factors as $(x - 1)\bigl(\bigl(y_0^2-1\bigr)xy - 2y_0y^2 +2y_0x+ \bigl(y_0^2-1\bigr)y\bigr)=0$,
\item[(3)] for $x_0= -1$, $y_0\ne 0$, \eqref{polartH} factors as $(x + 1)\bigl(\bigl(y_0^2-1\bigr)xy + 2y_0y^2 +2y_0x- \bigl(y_0^2-1\bigr)y\bigr)=0$.
\end{itemize}
The cases 2) and 3) give Type 11, the substitution $x\to -x$ reducing one to the other.

For infinite vertexes $p_0$, the pencil $x=k$ gives a non-smooth conic. Thus the pencil can be chosen as $y=ax+k$. A line from the pencil is tangent to the conic $c_r$ with
\smash{$
r=\frac{k^2+a^2+1}{2ak}$}.
 The points, where the lines $l_r$ touch $c_r$, are
\[
x=\frac{k\bigl(k^2 -a^2+1\bigr)}{a\bigl(k^2-a^2-1\bigr)}, \qquad y=\frac{2k}{a^2 - k^2 + 1}.
\]
Excluding $k$, we get the cubic
$
y^3+a^2x^2y+2axy^2+2ax+\bigl(1-a^2\bigr)y=0$.
The cubic splits into a line and a conic only for $a=0$ or for $a^2 \pm 2{\rm i}a -1=0$. The former case gives non-real and non-smooth conic $y^2+1=0$, the latter -- the non-real conic $xy \pm {\rm i} \bigl(y^2+2\bigr)=0$.

 If the hyperbolic line contains the point where the polar conic plane touches the Darboux quadric, then we sent this point to $(0,0,-1)$ and the polar line to $X=Y=0$. In the affine coordinates on the plane $z=-1$, the conics $c_r$ are concentric circles. The case of concentric circles was considered above. We get Types 12 and~13.

$\bullet$ {\it Polar conic plane tangent to Darboux quadric, elliptic pencil}.
If one vertex of the elliptic pencil coincides with the point where the polar conic plane touches the Darboux quadric, then the polar curve is planar.

Thus we can sent the tangent point to $(0,0,-1)$ and the vertexes to $(\pm 1,0,0)$. The conics $c_r$ have equations
$
\bigl(r^2 + 1\bigr)x^2 + (y+r)^2 = r^2+1$.
If the polar conic coincides with one of $c_r$, then the web curvature
\[
K_B=-\frac{64\bigl(r^2 + 1\bigr)^2x\bigl(y^2 + 1\bigr)\bigl( rx^2 + ry^2+ \bigl(3+r^2\bigr)y - r \bigr)\bigl(2ry - x^2 - y^2 + 1\bigr)^4}{\bigl(x^2 - y^2 - 1\bigr)^4\bigl(r^2 - x^2 + 1\bigr)^6}
\]
vanishes only for $r=\pm {\rm i}$, the conic $c_r$ being non-smooth for these values.

A line from the pencil $y=y_0+k(x-x_0)$ is tangent to the conic $c_r$ with
\[
r=\frac{\bigl(x_0^2 - 1\bigr)k^2 - 2x_0y_0k + \bigl(y_0^2 - 1\bigr)}{2(x_0k - y_0)}.
\]
 The points, where the lines $l_r$ touch $c_r$, are
\begin{gather*}
x=\frac{2k(x_0k - y_0)}{\bigl(x_0^2 + 1\bigr)k^2 - 2x_0y_0k + y_0^2 + 1},\\
y=-\frac{x_0\bigl(x_0^2 - 1\bigr)k^3 - y_0\bigl(3x_0^2 - 1\bigr)k^2 + x_0\bigl(3y_0^2 + 1\bigr)k - y_0\bigl(y_0^2 + 1\bigr)}{\bigl(x_0^2 + 1\bigr)k^2 - 2x_0y_0k + y_0^2 + 1}.
\end{gather*}
Excluding $k$, we get the cubic
\[
\bigl(y_0^2 + 1\bigr)x^3 - 2x_0y_0x^2y+ \bigl(x_0^2 + 1\bigr)xy^2 - 2x_0x^2 - 2x_0y^2 + \bigl(x_0^2 - y_0^2\bigr)x + 2x_0y_0y=0.
\]
This cubic splits into a conic and a line in 4
 cases:
 \begin{itemize}
\item[(1)] for $y_0=\pm {\rm i}$, the cubic equation factors as
\[
(\pm i-y)\bigl(\pm 2{\rm i}x_0x^2 -\bigl(1+x_0^2\bigr)xy \mp {\rm i}\bigl(1+x_0^2\bigr)x + 2x_0y\bigr)=0,
\]
\item[(2)] for $y_0=\pm {\rm i}x_0$, the cubic equation factors as
\[
 (x \pm {\rm i} y)\bigl(\bigl(1-x_0^2\bigr)x^2 \mp {\rm i} \bigl(x_0^2+1\bigr)xy - 2x_0x \pm2 {\rm i} x_0y + 2x_0^2\bigr)=0,
\]
\item[(3)] for $x_0=0$, the cubic equation factors as $x\bigl(\bigl(y_0^2 + 1\bigr)x^2 + y^2 - y_0^2\bigr)=0$,
\item[(4)] for $x_0=\pm 1$, the cubic equation factors as $(x\mp 1)\bigl(\bigl(y_0^2 + 1\bigr)x^2+2y^2 \mp 2y_0xy \pm \bigl(y_0^2 - 1\bigr)x - 2y_0y\bigr)=0$.
 \end{itemize}
In the cases 1) and 2) the conic is not real, the case 3) gives Type 14, for the case 4) the web curvature does not vanish.

For infinite vertexes $p_0$, the pencil $x=k$ gives a non-smooth conic. Thus the pencil can be chosen as $y=ax+k$. A line from the pencil is tangent to the conic $c_r$ with
\smash{$
r=\frac{a^2 - k^2 + 1}{2k}$}.
 The points, where the lines $l_r$ touch $c_r$, are
\[
x=-\frac{2ak}{k^2+a^2+1}, \qquad y=\frac{k\bigl( k^2-a^2 + 1\bigr)}{k^2+a^2+1}.
\]
Excluding $k$, we get the cubic
\[
a^2x^3 - 2ax^2y + xy^2 +\bigl(1-a^2\bigr)x + 2ay=0.
\]
The cubic splits into a line and a conic in 2 cases:
\begin{itemize}
\item[(1)] for $a=0$ the cubic equation factors as $x\bigl(y^2+1\bigr)=0$ and the conic is non-smooth
\item[(2)] for $a=\pm {\rm i}$ the cubic equation factors as $(x \pm {\rm i} y)\bigl(x^2 \pm {\rm i} xy-2\bigr)=0$, and the conic is not real.
\end{itemize}

$\bullet$ {\it Polar conic plane tangent to Darboux quadric, parabolic pencil}.
For non-planar polar curves, the points, where the polar line and polar conic plane touch the Darboux quadric, are different. Thus we can normalize the polar conic plane to $z=-1$ and the polar line to $x=0$, $z=1$. This configuration is preserved by $B_z$. The conics $c_r$ have equations
$
\frac{r}{4}x^2 + y =\frac{1}{r}$.
If the polar conic coincides with one of $c_r$, then the web curvature
\[
K_B=-\frac{4096r^5xy^2\bigl(ry + 2x^2 + 2y^2\bigr)\bigl(ry - x^2 - y^2\bigr)^4}{\bigl(x^2 - y^2\bigr)^4\bigl(r^2 - 4x^2\bigr)^6}
\]
vanishes only for $r=0$, the conic $c_r$ being non-smooth for this value.

 A line from the pencil $y=y_0+k(x-x_0)$ is tangent to the conic $c_r$ with
$
r=\frac{k^2+1}{y_0-x_0k}$.
 The points, where the lines $l_r$ touch $c_r$, are
\[
x=\frac{2k(x_0k - y_0)}{k^2 + 1},\qquad
y=\frac{x_0k^3 -y_0 k^2 -x_0 k + y_0}{k^2 + 1}.
\]
Excluding $k$, we get the cubic
\[
x^3+ xy^2 - 2x_0x^2- 2x_0y^2 + \bigl(x_0^2 - y_0^2\bigr)x + 2x_0y_0y=0.
\]
This cubic splits into a conic and a line in 2
 cases:
 \begin{itemize}
\item[(1)] for $y_0=\alpha x_0$, where $\alpha ^2 \pm 2{\rm i} \alpha -1=0$, the cubic equation factors as
\[
(x \pm {\rm i} y)\bigl(\pm {\rm i} xy -x^2 +2x_0x\mp 2{\rm i} x_0y - 2x_0^2\bigr)=0,
\]
\item[(2)] for $x_0=0$, the cubic equation factors as $x\bigl(x^2+y^2-y_0^2\bigr)=0$.
 \end{itemize}
 In the case 1), the conic is not real, the case 2) gives Type 15 after rescaling by $B_z$.

 For infinite vertexes $p_0$, the pencil $y=k$ gives a non-smooth conic. Thus the pencil can be chosen as $x=ay+k$. A line from the pencil is tangent to the conic $c_r$ with
$
r=-\frac{a^2+ 1}{ak}$.
 The points, where the lines $l_r$ touch $c_r$, are
\[
x=\frac{2k}{a^2+1}, \qquad y=\frac{k\bigl(1- a^2 \bigr)}{a\bigl(1+a^2\bigr)}.
\]
Excluding $k$, we get the line
$\bigl(a^2 - 1\bigr)x + 2y=0$.

 Thus we have considered all the cases with non-planar polar curve, the theorem is proved.
\end{proof}

\begin{Corollary}
Suppose that the polar curve of a hexagonal circular $3$-web is reducible algebraic of degree three. Then the polar curve is either planar or the web is M\"obius equivalent to one described by Theorems {\rm\ref{3pencils}} and {\rm\ref{conicPencil}}.
\end{Corollary}

\section{Hexagonal circular 3-webs with M\"obius symmetries}\label{symetric webs}
 A M\"obius transformation is a symmetry of a circular 3-web if it maps any web circle to a web circle.
The M\"obius group in $\mathbb{RP}^3$ can be realized as $PGL_2(\mathbb{C})$, or equivalently, as the group of fractional-linear transformations of $z=x+{\rm i}y$, where $(x,y)$ are cartesian coordinates in the plane. A generator of any 1-dimensional subalgebra can be brought to the Jordan normal form by adjoint action. The generator can be chosen either as
$
\big(
\begin{smallmatrix}
0&1\\
0&0
\end{smallmatrix}
\big)$, or $
\big(
\begin{smallmatrix}
\lambda &0\\
0&-\lambda
\end{smallmatrix}
\big)$,
where $\lambda=\alpha+{\rm i}\beta$ is some complex number with $\operatorname{Re}(\lambda)=\alpha$, $\operatorname{Im}(\lambda)=\beta$.

One way to obtain symmetric hexagonal 3-webs is provided by the Wunderlich construction (see \cite{W-38} and Introduction).

Another easy way to produce hexagonal circular 3-webs is to choose two orbits of polar points such that one is a conic and the other is a coplanar straight line, these two orbits forming a~polar curve. This construction may degenerate if there are 3 orbits which are coplanar straight lines.

For translations $(x,y)\mapsto (x,y+u)$, the orbit of a polar point for a circle $(x-a)^2+y^2=r$, parametrized by $u$ as follows $\bigl[-2a:-2u: a^2+u^2-r-1: -a^2-u^2+r-1\bigr]$, is a conic in the plane $-X+a(Z+U)=0$. For a nonvertical line $y=ax$, the orbit is a line $Z+U=X+aY=0$, parametrized by $u$ via $[a:-1:-u:u]$.
Thus we immediately get the following hexagonal 3-webs symmetric by translations.

 {\bf T1.\ 3~families of parallel lines.}
 M\"obius orbits of such webs form a two-parametric family. The polar curve splits into 3 coplanar lines. Observe that any web of this family has a~3-dimensional symmetry group.

 {\bf T2.\ Polar curve splits into conic and coplanar line.}
 There is only one M\"obius class of such webs, any representative is formed by horizontal lines $y=\text{const}$ and by the orbit of a~circle which can be chosen as the unitary one centered at the origin.

 {\bf T3.\ Wunderlich's type.}
 There are several types, depending on the position of generating curves.
\begin{enumerate}\itemsep=0pt
\item[(1)] Nondegenerate type. Webs are formed by vertical lines and by two different orbits of circles. The family of orbits is two-parametric: by translation and rescaling we can fix one orbit.

\item[(2)] Coinciding circle orbits. There is only one M\"obius class of this type. It was already obtained as Type 9 in the classification of Theorem \ref{conicPencil}.

\item[(3)] Generated by a circle and a line. The family of orbits is one-parametric: by translation and rescaling we can fix the generating circle.
\end{enumerate}

 For dilatations $(x,y)\mapsto (ux,uy)$, the orbit of a polar point for a circle $(x-a)^2+(y-b)^2=r$, parametrized by $u$ via $\bigl[-2au:-2bu: a^2+b^2-r-u^2: -\bigl(a^2+b^2-r\bigr)-u^2\bigr]$, is a conic in the plane $bX-aY=0$ if $(a,b)\ne (0,0)$ and $a^2+b^2\ne r$, the hyperbolic line $X=Y=0$ if $(a,b)= (0,0)$, and the parabolic line $bX-aY=U-Z=0$ if $a^2+b^2=r$. For a line $ax+by+c=0$ with $c\ne 0$, the orbit is the parabolic line $bX-aU=Z+U=0$, parametrized by $u$ via $[au:bu:c:-c]$. The lines with $c=0$ are invariant.
Invoking the classification of the webs with 3 pencils, we list the following hexagonal 3-webs symmetric by dilatations.

 {\bf D1=T1.\ 3 families of parallel lines.}

 {\bf D2.\ 2 coplanar parabolic lines and hyperbolic line intersecting them.}
 Family of parallel lines, the parabolic pencil with circles tangent to a line of the family through the origin, and the family of concentric circles with the center at the origin. There is only one M\"obius class of such webs.

 {\bf D3.\ 2 dual parabolic lines and hyperbolic line through their common point.}
 2~orthogonal families of parallel lines and the family of concentric circles with the center at the origin. There is only one M\"obius class of such webs, we have already presented it on the right of Figure~\ref{f2}.

 {\bf D4.\ Polar curve splits into conic and coplanar hyperbolic line.}
 Family of circles obtained by the dilatations from one not passing through the origin and not having its center at the origin and the family of concentric circles with the center at the origin. M\"obius orbits of such webs form a one-parameter family.

 {\bf D5.\ Polar curve splits into conic and tangent parabolic line.}
 Family of circles obtained by dilatations from one not passing through the origin and not having its center at the origin and the family of parallel lines orthogonal to the orbit of centers of the circles. M\"obius orbits of such webs form a one-parameter family.

 {\bf D6.\ Wunderlich's type.}
\begin{enumerate}\itemsep=0pt
\item[(1)] Nondegenerate type, polar curve with 2 conics and elliptic line. Webs are formed by pencil of lines centered at the origin and by two different orbits of circles at general position. The family of M\"obius orbits is 3-parametric: one can choose the centers of circles on a fixed circle centered at the origin and normalize by rotations.

\item[(2)] Polar curve with conic and elliptic line. Webs are formed by pencil of lines centered at the origin and orbit of a circle. The family of M\"obius orbits is 1-parametric, it is Type 5 in the classification of Theorem \ref{conicPencil}.

\item[(3)] Polar curve with conic, hyperbolic and elliptic line. Webs are formed by pencil of lines centered at the origin, orbit of a circle, and the family of concentric circles with the center at the origin. The family of M\"obius orbits is 1-parametric.

\item[(4)] Polar curve with conic, parabolic and elliptic line. Webs are formed by pencil of lines centered at the origin, orbit of a circle, and a family of parallel lines. The family of M\"obius orbits is 2-parametric.

\item[(5)] Polar curve with hyperbolic line, elliptic line dual to hyperbolic, and parabolic line intersecting them. Webs are formed by the pencil of lines centered at the origin, the family of concentric circles with the center at the origin, family of parallel lines. There is only one M\"obius type. This web is shown in the center of Figure~\ref{f2}.

\item[(6)] Polar curve with 2 parabolic lines touching Darboux quadric at the same point and coplanar elliptic line. Webs are formed by the pencil of lines centered at the origin and 2 families of parallel lines. The family of M\"obius orbits is 1-parametric.

\item[(7)] Polar curve with 2 parabolic lines and elliptic line whose dual joins the touching points of parabolic lines with Darboux quadric. Webs are formed by the pencil of lines centered at the origin, a family of parallel lines and a parabolic pencil with the vertex at the origin. The family of M\"obius orbits is 1-parametric. A representative of this web is shown in Figure~\ref{f4}.

\end{enumerate}

 For rotations $(x,y)\mapsto (x\cos(t)+y\sin(t),-x\sin(t)+y\cos(t))$, the orbit of a polar point for a~circle $(x-a)^2+y^2=r$, parameterized by $t$ via $\bigl[-2a\cos(t):-2a\sin(t): a^2-r-1: -a^2+r-1\bigr]$, is a circle in the plane $\bigl(a^2-r+1\bigr)Z+\bigl(a^2-r-1\bigr)U=0$ if $a\ne 0$. For a line $ax+by+c=0$ with~${c\ne 0}$, the orbit of its polar is also a circle $[a\cos(t)-b\sin(t):a\sin(t)+b\cos(t):c:-c]$ in the plane $Z+U=0$. The circles with $a=0$ are invariant, the orbit of a polar point of a line with $c=0$ is an elliptic line.
Thus one can easily list the following hexagonal 3-webs symmetric by rotations.

 {\bf R1.\ Polar curve is a circle symmetric by rotations around $\boldsymbol{ z}$-axis and coplanar elliptic line.}
 The family of M\"obius orbits is 1-parametric.

 {\bf R2.\ Wunderlich's type.}
 One foliation is formed by orbits of points by rotations and the other two are images of two circles by rotation, these circles can coincide or ``degenerate'' into straight lines. Examples of such webs with coinciding generators are shown in Figure~\ref{123} (center), Figure~\ref{456}~(left), and Figure~\ref{131415}~(left). The reader easily describes different types and computes corresponding M\"obius orbit dimension of the webs.

\begin{Theorem}
Hexagonal circular $3$-web with $1$-dimensional M\"obius symmetry is M\"obius equivalent to one of the above described T-, D-, or R-types.
\end{Theorem}
\begin{proof}
 A circular 3-web, symmetric by a given 1-parametric subgroup of the M\"obius group and not obtained by the Wunderlich construction, is fixed by a choice of three curves, each being either a circle or a straight line, one from each foliation. Two circles can coincide: an orbit of circle still gives a (singular) 2-web. Moreover, one can move around these curves by the stabilizer of the infinitesimal generator of the subgroup.

$\bullet$ {\it Webs symmetric by translations $\partial_y$}.
Consider a point such that none of the leaf tangents is parallel to the field $\partial_y$. If a generating curve $C_1$ of some foliation is a circle, then there are 2 circles from the orbit of $C_1$ passing through this point. Thus there are two locally defined direction fields $\partial_x+P_{\pm}\partial_y$ tangent to these 2 circles. Globally they are not separable: one direction swaps for the other upon running along $C_1$.

Consider one of the foliations and the corresponding direction field $\partial_x+P\partial_y$. Since the foliation is symmetric, the slope $P$ does not depend on $y$. Since all the leaves are circles of the same curvature, the foliation has a first integral
\smash{$
\frac{({P'})^2}{(P^2+1)^3}$}, where $P'=\frac{{\rm d}P}{{\rm d}x}$.
The differential equation
\smash{$
\frac{({P'})^2}{(P^2+1)^3}=A^2=\text{const}
$}
has general solution
\[
P(x)=\frac{A(x-x_0)}{\sqrt{1-A^2(x-x_0)^2}},
\]
the zero value of $A$ corresponding to the foliation by straight parallel lines.
Let the slopes of the other two web foliations be $Q$ and $R$, and the connection form be
$\gamma=\alpha(x){\rm d}x+\beta(x){\rm d}y$. Hexagonality condition ${\rm d}\gamma=0$ implies $\beta(x)=\text{const}$, which amounts to
\begin{equation}\label{transconect}
k=\frac{P'}{(P-Q)(P-R)}+\frac{Q'}{(Q-P)(Q-R)}+\frac{R'}{(R-P)(R-Q)}=\text{const}.
\end{equation}
If all foliations are formed by circles and not by straight lines, then the above identity is not possible. In fact, each of the slopes $P$, $Q$, $R$, considered as function of complex $x$, has two ramification points, for example, $P$ has singularities at $x_{\pm}=x_0\pm 1/A$. If at least one of the 6~ramification points is not coinciding with one of the others, then, supposing it be $x_{+}$ of $P$ and expanding the expression~\eqref{transconect} for $x=x_{+}+t^2$ by $t$ at $t=0$, one sees that it has a simple pole and therefore cannot be constant. Therefore, each ramification point of any slope coincides with a ramification point of another slope.

The ramification points correspond to the group orbits (lines $x=\text{const}$) tangent to the foliation circles. If only two circles are tangent along some orbit $x=\text{const}$, then, applying Lemma~\ref{sing}, we conclude that either this orbit belongs to the third foliation and the web is of Wunderlich's type or the two generating circles coincide. In the latter case, the ramification points of the third foliation again must coincide with the common ramification points of the first two. This means that all orbits of generating circles coincide and we do not have a 3-web.

Therefore, either all leaves of all foliations are straight lines or one of the foliation, say the one corresponding to $R$, is formed by straight lines with slope $R=\text{const}$ and the other two are formed by orbits of the same circle and $Q=-P$. Then it is immediate that $k=0$ and $R=0$. We get the type $T2$. The webs of Wunderlich's type, which are hexagonal, are excluded from consideration by the coordinate choice.

$\bullet$ {\it Webs symmetric by dilatations $x\partial_x+y\partial_y$}.
Stabilizer of the 1-dimensional algebra spanned by $x\partial_x+y\partial_y$ is generated by rotation $y\partial_x-x\partial_y$ and dilatation $x\partial_x+y\partial_y$. The orbit of the polar of the generating curve is either a conic, or a hyperbolic line, or a parabolic line. Parabolic lines touch the Darboux quadric at one of the stationary point of the dilatation.

We use modified polar coordinates rectifying the symmetry
\begin{equation}\label{polar}
u = \arctan\left(\frac{y}{x}\right), \qquad v= \frac{1}{2} \ln\bigl(x^2 + y^2\bigr).
\end{equation}
A curve $x\mapsto (x,y(x))$ is a circle if and only if
\smash{$
y'''=\frac{3y'(y'')^2}{1+(y')^2}$},
therefore integral curves of a~symmetric vector field $\partial_u+P(u)\partial_v$ are circles
if and only if
\begin{equation}\label{scale2order}
P''=\frac{3P}{P^2 + 1}(P')^2 - P\bigl(P^2 + 1\bigr).
\end{equation}

The integral curves are straight lines if and only if $y''=0$, which is equivalent to
$P'=P^2+1$. Hence \big(apply $z\mapsto \frac{1}{z}$\big) the integral curves are circles of a parabolic pencil with the vertex at the origin if and only if $P'=-\big(P^2+1\big)$ thus giving $P(u)=\tan(u-u_0)$ and $P(u)=-\tan(u-u_0)$ respectively.
The second-order equation \eqref{scale2order} has a first integral
\[
\frac{(P')^2}{\bigl(P^2 + 1\bigr)^3} - \frac{1}{P^2 + 1} =A=\text{const},
\]
allowing also to integrate it
\[
P(u)=\frac{\sqrt{A + 1}\tan(u - u_0)}{\sqrt{1 - A\tan^2(u - u_0)}}.
\]
To the hyperbolic pencil with the circles centered at the origin corresponds the solution $P\equiv 0$.

 For hexagonal webs, the connection form is
$\gamma=\alpha(u){\rm d}u+\beta(u){\rm d}v$. Hexagonality condition~${{\rm d}\gamma=0}$ implies $\beta(u)=\text{const}$, giving again \eqref{transconect} (where the slopes of the other foliations are~$Q$ and $R$).
Observe that the choice of local coordinates $u$, $v$ excludes webs of Wunderlich's type and we have to show that hexagonal are only the types D1, D2, D3, D4, D5.
This can be done as follows. The types D1, D2, D3 are webs with three pencils, the case being settled earlier. Therefore, we study the webs whose polar curve includes at least one conic.

The approach used for translation works also in this case if we pass to complex webs: the singular points of solutions $P$, $Q$, $R$ become complex if the corresponding generating circle on the Darboux quadric separates the stable points $(0,0,\pm 1)$ of the dilatation.
Considering behavior of the expression $k$ at singular points of $P$, $Q$, $R$, we see that necessarily a singular point of slope for one foliation must coincide with a singular point for another. Singular points are the points (possibly complex) where the group orbit touches the generating circle, or lies either on a line of one foliation or on the common tangent to circles of parabolic pencil with the vertex at $(0,0)$.

Applying Lemma \ref{sing} as in the case of translation, we infer that, for webs of non-Wunderlich's type, having two foliations with conic polar curves and coinciding singular points, these polar curves must coincide and the third polar must be a line.

Consider the points where the circles from the common orbit are tangent. Lemma \ref{sing} implies that for non-Wunderlich's type the leaves of the third foliations must be also tangent to the circles and we get either type D4 or D5.

Finally, if two components of the polar curve are lines and the third is a conic then the lines must be parabolic for
non-Wunderlich's type. In fact, considering the points where circles of the hyperbolic pencil are tangent to circles corresponding to conic, we conclude by Lemma \ref{sing} that the third foliation lines are also tangent to the circles at these points but then the singular points of the ``conic'' solution to \eqref{scale2order} can not be compensated.

Thus the polar lines are parabolic. The expression $k$ cannot be constant if the poles of $Q$, $R$, corresponding to these lines, do not coincide with singular points of $P$, which represent a conic in the polar curve. Therefore, the generating circle of this ``conic'' solution does not separate the stationary point of the dilatation and the singular points are real. Then by Lemma \ref{sing} the web cannot be hexagonal.

$\bullet$ {\it Webs symmetric by rotations $y\partial_x-x\partial_y$}.
Stabilizer of the 1-dimensional subalgebra spanned by $y\partial_x-x\partial_y$ is generated by rotations $y\partial_x-x\partial_y$ and dilatations $x\partial_x+y\partial_y$. The polar orbit for a generating curve is either a circle or an elliptic line.

We again use polar coordinates \eqref{polar} to describe symmetric vector fields $\partial_v+P(v)\partial_u$.
The integral curves of such vector field are circles (in coordinates $x$, $y$, of course)
if and only if
\begin{equation}\label{rot2order}
P''=\frac{3P}{P^2 + 1}(P')^2 + P\bigl(P^2 + 1\bigr).
\end{equation}
The integral curves are straight lines if and only if
$P'=-P\bigl(P^2+1\bigr)$. The integral curves are circles passing through the origin if and only if $P'=P\bigl(P^2+1\bigr)$. Solutions of these two equations are \smash{$P(v)=\frac{1}{\sqrt{e^{2(v-v_0)}-1}}$} and \smash{$P(v)=\frac{1}{\sqrt{e^{-2(v-v_0)}-1}}$}, respectively.
The second-order equation \eqref{rot2order} has a first integral
\[
\frac{(P')^2}{\bigl(P^2 + 1\bigr)^3} + \frac{1}{P^2 + 1} =A^2=\text{const},
\]
allowing also to integrate it
\[
P(v)=\frac{\sqrt{A^2 - 1}\tanh(v - v_0)}{\sqrt{1 - A^2\tanh^2(v - v_0)}}, \qquad A^2\ne 1.
\]
The value $A^2=1$ gives the special solutions with generating circles passing through the stationary points of the symmetry, i.e., elliptic pencil with the lines passing through the origin with~${P\equiv 0}$.
Analysis of the behavior of $k$ at singular points and use of Lemma \ref{sing}, similar to the ones performed above, show that only the type R1 is hexagonal, the Wunderlich types being excluded by the choice of variables. (In fact, multiplying by $i$ the independent variable of the differential equation \eqref{scale2order} reduces it to \eqref{rot2order}.)

$\bullet$ {\it Loxodromic symmetry.}
Finally, we show that there is no hexagonal 3-webs symmetric by loxodromic vector field $y\partial_x-x\partial_y+\kappa (x\partial_x+y\partial_y)$ for any $\kappa\ne 0$. Note that the Wunderlich construction does not give circular webs as the symmetry orbits are spirals. We use the following coordinates:
\[
s= \frac{\kappa}{2} \ln\bigl(x^2 + y^2\bigr)-\arctan\left(\frac{y}{x}\right), \qquad t=\kappa\arctan\left(\frac{y}{x}\right)+ \frac{1}{2} \ln\bigl(x^2 + y^2\bigr),
\]
the variable $t$ being invariant by the symmetry.
Integral curves of a symmetric vector field $\partial_t+P(t)\partial_s$ are circles (in coordinates $x$, $y$)
if and only if
\begin{equation}\label{drom2order}
P''=\frac{3P}{P^2 + 1}(P')^2 + \frac{\bigl(P^2 + 1\bigr)(\kappa P + 1)(P-\kappa)}{\bigl(\kappa^2 + 1\bigr)^2}.
\end{equation}
This equation has a first integral
\[
\frac{(P')^2}{\bigl(P^2 + 1\bigr)^3} + \frac{2\kappa P -\kappa^2+1}{\bigl(\kappa ^2+1\bigr)\bigl(P^2 + 1\bigr)} =A=\text{const}.
\]
This integral does not allow to integrate \eqref{drom2order} in elementary functions but allows to study the behavior of solutions at singular points. A singular point emerges when a symmetry orbit is tangent to a circle (or a line) of corresponding foliation.

Consider possible singularity types. If we exclude from consideration the stationary points of symmetry, then the symmetry vector field touches a generic circle at two points and the tangency is simple. The corresponding solution $P$ has two singularities if the tangency points belong to different orbits and only one if the points lie on the same orbit. There are circles, for which two tangency points merge to give only one singularity. If the generating curve is a line or a circle through the origin, then the solution has only one singularity. Finally, there are two constant solutions, namely $P=-1/\kappa$ corresponding to the invariant hyperbolic line $X=Y=0$ and~${P=k}$ corresponding to the invariant elliptic line $Z=U=0$.

If a solution $P$ has two singularities at $t_1$, $t_2$, then $A\ne 0$ and the singularity is of the same type as for the non-loxodromic cases
\[
P(t)=\frac{1}{\sqrt[4]{4A}\sqrt{t-t_i}}+\bigl\{\text{analytic function of } \sqrt{t-t_i}\bigr\}.
\]

If a solution $P$ is generated by a circle tangent to the symmetry trajectory $t=t_0$ and the tangency is of second order, then there is only one singularity of the following type at~$t=t_0$ \smash{$
P(t)=c(t-t_0)^{-\frac{2}{3}}+\cdots$},
where $c\ne 0$ and the omitted terms are not essential for our analysis.
The condition of double tangency is equivalent to $A=0$. The corresponding generating circle $(x-a)^2+y^2=r^2$ verifies the relation \smash{$r^2=\frac{(k^2 + 1)a^2}{k^2}$}. For an orbit $t=t_0$, there is at most one such circle.

Let $P$, $Q$, $R$ be solutions giving a hexagonal 3-web. These solutions, as well as the coordinates~$s$,~$t$, are defined only locally but we can prolong them along a symmetry orbit, along a leaf of some of the 3 foliations, or along any curve, as long as we do not meet a singular point of one of $P$, $Q$, $R$.
The condition of hexagonality \eqref{transconect} remains satisfied along any such prolongations. Therefore, we cannot meet a singularity of only one of $P$, $Q$, $R$, they emerge necessarily at least in pairs.

Suppose there is a symmetric hexagonal 3-web. Consider a non-singular point. There are 3 leaves passing through it. Each leaf can be considered as the generating curve of the respective foliation. At least one of the corresponding solutions $P$, $Q$, $R$ has a singular point. Let us run along the respective leaf until we meet a singularity $t_1$ of $P$, $Q$ or $R$. Then at least two of~$P$,~$Q$,~$R$ are singular at $t_1$ and the orbit $t=t_1$ is tangent to at least two web leaves at a singular point~$p_s$. Since the symmetry orbit is not a circle, Lemma \ref{sing} implies that either exactly two leaves at $p_s$ are tangent and therefore coincide or all three leaves are tangent at $p_s$.

In the former case, suppose that the coinciding leaves are $C_Q$ and $C_R$. Then the leaf $C_P$ at~$p_s$ is different from $C_Q=C_R$ at $p_s$. Let us go along $C_P$ keeping track of $Q$, $R$ until one of the two leaves corresponding to $Q$ and $R$ touches $C_P$ at some $\bar{p}_s$. Such point exists until all the foliations are formed by straight lines. Then by Lemma~\ref{sing} this leaf coincide with $C_P$ at $\bar{p}_s$, thus all three generating curves of the web coincide and there is no 3-web. If all 3 generating curves are straight lines, then $C_P$ is necessary the line through the origin and $P=\kappa$. Applying the map~${z\mapsto 1/z}$, we transform the lines $C_Q$ and $C_R$ into circles and the above argument applies.

If all three leaves are tangent at $p_s$, then at least one leaf, say $C_P$, is different from any of the other two. Let us go again along $C_P$ keeping track of $Q$, $R$ until one of the two leaves corresponding to $Q$ and $R$ touches $C_P$. Let it be $C_Q$. Repeating the above used argument we conclude that such point $\bar{p}_s$ exists and $C_P$ coincides at $\bar{p}_s$ with~$C_Q$. Then $C_P$ coincides with~$C_Q$ also at $p_s$. Now either all 3 generating curves $C_P$, $C_Q$, $C_R$ coincide at $p_s$ and we do not have 3-web, or $C_R$ is different from $C_P=C_Q$ at $p_s$. Now we repeat the trick with prolongation, this time along $C_R$, and conclude that $C_P=C_Q=C_R$ at $p_s$. Thus all three generating curves coincide and we can get at most 2-web.
\end{proof}

 \section{Concluding remarks}

\subsection{Circular hexagonal 3-webs on surfaces}

Pottmann, Shi, and Skopenkov \cite{PSS-12} classified circular hexagonal 3-webs on nontrivial Darboux cyclides: such surfaces carry up to 6 one-parameter families of circles, 3 families can be picked up in 5 different ways to form a hexagonal web. In fact, as was proven by Lubbes \cite{L-21}, if through a general point of a surface in $\mathbb{R}^3$ pass at least 3 circles then the surface is either a plane, or a~sphere, or a Darboux cyclid.

\subsection[Erdogan's approach to Theorem 4.5]{Erdo\v{g}an's approach to Theorem \ref{3pencils}}

The first attempt to prove Theorem \ref{3pencils} appeared in \cite{E-89}. The idea was to choose a M\"obius normalization sending one of the vertexes to infinity, to set $y=0$, and to obtain ``sufficient number of equations'' to fix the pencil configurations. The author claimed that the curvature equation for $y=0$ is a polynomial one of degree 6 in $x$, though the calculation itself was not present. Nowadays, armed with a powerful computer (32GB of RAM is enough) and a symbolic computation system like Maple, one can perform this computations and check that the degree may be much higher (in fact, up to 18) for some choices of polar line types.

Anyway, brute computer force does work: with the above mentioned equipment the author of this paper managed to derive the classification results.
The treating has the following steps:
\begin{enumerate}\itemsep=0pt
\item[(1)] choosing an initial M\"obius normalization,
\item[(2)] computing the curvature,
\item[(3)] isolating and factoring the highest homogeneous part of the curvature equation,
\item[(4)] M\"obius renormalization adjusted to the geometric information obtained in the previous step and repeating from the step 2 until one makes the curvature vanish.
\end{enumerate}

\subsection{Boundaries of regular domain for hexagonal 3-webs}

To avoid heavy computation of the curvature in proving Theorem \ref{3pencils}, the author of \cite{S-05} suggested to use the structure of web singular set (see Lemma~\ref{sing}). The presented proof was not correct. The author argued that the web equation $u_3=F(u_1,u_2)$, relating first integrals $u_i$ of the web foliations $\mathcal{F}_i$, may be rewritten as $u_3=f(\alpha(u_1)+\beta(u_2))$ and used this form on the curve of singular points $\Gamma_1$. This argument is definitely wrong as the functions $\alpha$, $\beta$ typically have singularities on $\Gamma_1$: a simple counterexample is the web equation $u_3=u_1u_2$.

\subsection{Hexagonal 3-subwebs}

Consider the autodual tetrahedron with vertexes at $[1:0:0:0]$, $[0:1:0:0]$, $[0:0:1:0]$ and $[0:0:0:1]$. Lines joining the vertexes of this tetrahedron give a 6-web $A_6$ with 6 pencils of circles. Any 3-subweb of this 6-web is hexagonal: any 3 lines are either coplanar or give a~polar curve of a hexagonal 3-web. By direct computation (better use computer!) one checks the following claim.
\begin{Proposition}
The rank of the $6$-web $A_6$ is maximal, i.e., is equal to $10$.
\end{Proposition}
Another remarkable feature of this autodual 6-web is that the infinitesimal operators, corresponding to the pencils in the sense of Proposition \ref{operline}, form the basis $R_x$, $R_y$, $R_z$, $B_x$, $B_y$, $B_z$ of the Lie algebra of the M\"obius group so that the commutator of any two of them is either zero or an operator of the basis.

There is also autodual 4-web $A_4$ whose polar curve is the union of 4 pencils corresponding to $R_z$, $B_z$, $B_x-R_y$, $R_x+B_y$. This web has similar properties: any its 3-subweb is hexagonal, the commutator of any two of the four operators is either zero or an operator of the set, its rank is maximal.
One finds more examples with hexagonal subwebs among symmetric webs of Wunderlich's type (see Section \ref{symetric webs}).

\subsection{Conjecture}
The polar curve of a hexagonal circular 3-web is an algebraic curve such that
each its irreducible component is either a twisted cubic or a planar curve of degree at most 3. The author thanks an anonymous referee for providing an example of hexagonal circular 3-web with polar twisted cubic.

\appendix
\section{Appendix: Curvature equation}\label{AppendixA}
\begin{gather*}
\bigl(A^2 - 4B\bigr)\bigl(R^2 +AR+ B\bigr)\bigl[(AR + 2B)A_{xx}+A\bigl(R^2-B\bigr)A_{xy}-BR(A + 2R)A_{yy}
\\
\qquad{}-(A + 2R)B_{xx}+ \bigl(A^2-2R^2- 2B\bigr)B_{xy}+R\bigl(A^2 -2B+AR\bigr)B_{yy}
\\
\qquad{}+\bigl(4B-A^2\bigr)R_{xx}+A\bigl(A^2 - 4B\bigr)R_{xy}+B\bigl(4B-A^2\bigr)R_{yy}\bigr]
\\
\qquad{}+\bigl(A^2 - 4B\bigr)^2(A + 2R)R_x^2+(A + 2R)\bigl(A^2 - 4AR - 4R^2 - 8B\bigr)B_x^2+
\\
\qquad{} - \bigl(2A^3R^2 + A^2R^3 + 7A^2BR + 4ABR^2 + 4BR^3 + 4AB^2 - 4B^2R\bigr)A_x^2
\\
\qquad{}+\bigl(B-R^2\bigr)\bigl(2A^3R + A^2R^2 + A^2B + 4BR^2 + 4B^2\bigr)A_xA_y \\
\qquad{}- A\bigl(A^2 - 4B\bigr)^2(A + 2R)R_xR_y
\\
\qquad{}+\bigl( 2A^3R-A^4 + 4A^2R^2 - 8AR^3 - 8R^4 + 8A^2B - 16BR^2 - 8B^2\bigr)B_xB_y
 \\
\qquad{}+BR\bigl(2A^3R + 7A^2R^2 + 4AR^3 + A^2B + 4ABR - 4BR^2 + 4B^2\bigr)A_y^2+
 \\
\qquad{}-R\bigl(A^4 - A^3R - 6A^2R^2 - 4AR^3 - 8A^2B - 8ABR + 8B^2\bigr)B_y^2
 \\
\qquad{}+B\bigl(A^2 - 4B\bigr)^2(A + 2R)R_y^2 +\bigl(A^2 - 4B\bigr)\bigl(A^2R - AR^2 - AB - 8BR\bigr)A_xR_x
 \\
\qquad{}+\bigl(3A^3R + 13A^2R^2 + 8AR^3 + A^2B + 12ABR - 4BR^2 + 12B^2\bigr)A_xB_x
 \\
\qquad{}+2\bigl(A^2 - 4B\bigr)\bigl(A^2 + AR + R^2 - 3B\bigr)B_xR_x
+\bigl( 5A^2R^3-A^4R - A^3R^2 + 4AR^4 \\
\qquad{} + A^2BR+ 4ABR^2 - 4BR^3 - 4AB^2 - 4B^2R\bigr)[A_xB_y+A_yB_x]
 \\
\qquad{}+\bigl(A^2 - 4B\bigr)\bigl(A^2R^2 + A^2B + 2ABR - 2BR^2 - 2B^2\bigr)[A_xR_y+A_yR_x] +
 \\
\qquad{}-A\bigl(A^2 - 4B\bigr)\bigl(A^2 + AR + R^2 - 3B\bigr)[B_xR_y+B_yR_x]
 \\
\qquad{}+B\bigl(A^2 - 4B\bigr)\bigl(A^2R - AR^2 - AB - 8BR\bigr)A_yR_y
 \\
\qquad{}+\bigl(4B-A^2 \bigr)\bigl(A^3R + A^2R^2 - A^2B - 6ABR - 6BR^2 + 2B^2\bigr)B_yR_y+
 \\
\qquad{}-R\bigl(2A^4R + 5A^3R^2 + 3A^2R^3 - A^2BR + 4ABR^2 + 4BR^3 + 8AB^2 + 20B^2R\bigr)A_yB_y \\
\qquad{}=0.
\end{gather*}

\subsection*{Acknowledgements}
This research was supported by FAPESP grant \# 2022/12813-5.
The author thanks the anonymous referees for valuable suggestions.

\pdfbookmark[1]{References}{ref}
\LastPageEnding

\end{document}